\documentclass[12pt]{article}

\usepackage[utf8]{inputenc}
\usepackage[english]{babel}

\usepackage{ifpdf}
\usepackage{latexsym}
\usepackage{amssymb}
\usepackage{pifont}
\usepackage{amsopn}
\usepackage{amsthm}
\usepackage{bbold}
\usepackage{graphicx}
\usepackage{tikz}
\usepackage[font={small,it}]{caption}
\usepackage{float}

\usepackage[pdftex,pdfstartview=FitH,pdfborderstyle={/S/B/W 1}, colorlinks=true,linkcolor=blue,urlcolor=blue,citecolor=blue,pagebackref=false]{hyperref}

\usepackage{bookmark}

\usepackage[capitalize,nameinlink]{cleveref}

\newtheorem{theorem}{Theorem}[section]
\newtheorem{proposition}[theorem]{Proposition}
\newtheorem{lemma}[theorem]{Lemma}
\newtheorem{corollary}[theorem]{Corollary}
\newtheorem{definition}[theorem]{Definition}
\newtheorem{notation}[theorem]{Notation}

\newcommand{\rr}{\mathbb{R}}
\newcommand{\cc}{\mathbb{C}}

\newcommand{\psltwo}{\mathrm{PSL}_2}

\begin{document}

\pagestyle{plain}

\title{CAT(-1)-Type Properties for Teichm\"uller Space}
\author{Ian Frankel}

\maketitle

\begin{abstract}We show that if Teichm\"uller geodesics spend enough time in the thick part of moduli space, they display CAT(-1)-type properties. In particular, they exponentially contract along strongly stable leaves. As an application we prove two closing lemmas.\end{abstract}

\tableofcontents

\section{Introduction and Statement of Main Results}

\noindent A \emph{half-translation surface} is an oriented surface $S$ with a maximal collection of charts to $\cc$ covering all but finitely many points of $S$ and whose transition maps are all of the form $z \mapsto z + C$ or $z \mapsto z - C$, where $C$ is a constant depending on the transition map. Such a surface naturally equipped with a Riemannian metric of zero curvature from the charts, which in general does not extend smoothly to the omitted points. A \emph{translation surface} is defined similarly, except the transition maps are all translations.\\

\noindent In case the underlying surface is compact and the Riemannian metric has finite area, the metric has a singular extension to the points not covered by the system of charts with cone-type singularities whose cone angles are all positive integer multiples of $\pi$. A finite area half-translation surface has \emph{type} $g,n$ if the underlying surface has genus $g$ and a collection of $n$ marked points, which includes every cone point of cone angle $\pi$. The remaining cone points may or may not be cone points for the metric. Two such surfaces $X_1,X_2$ are considered to be the same if their is an isometry between them, preserving the sets of marked points, and with the additional property the isometry takes the system of charts defining the half-translation structure for $X_1$ to the system of charts for the half-translation structure for $X_2$, or equivalently, if the $\cc$-derivatives of transition maps are constant and equal to $\pm 1$ for each pair of charts.\\

\noindent We refer to marked points and cone points as \emph{singularities}.\\

\noindent An example of such a surface is given by a \emph{pillowcase}: Consider a Euclidean rectangle and its reflection about one of its sides. The disjoint union of these two rectangles with each side identified with its image under reflection is a half-translation surface of type $0,4$. The four vertices of the original rectangle become cone points of angle $\pi$.\\

\noindent Every half-translation surface admits a cell structure whose 0-cells are the set of singularities, whose open 1-cells are geodesics not containing any singularity, and whose open 2-cells are open triangles with respect to charts defining the half-translation surface. Such a cell structure will henceforth be called a \emph{triangulation}. The 1-cells will be called \emph{saddle connections}.\\

\noindent The collection of all half-translation surfaces of type $g,n$ fit into a moduli space, which is endowed with various structures (algebro-geometric, complex analytic, measure-theoretic). This space will be called $QD(\mathcal{M}_{g,n})$ With respect to these structures, the typical point has $4g-4 + n$ cone points of cone angle $3\pi$ and $n$ cone points of cone angle $\pi$. At such a point, the saddle connections forming the $1$-skeleton can be assigned complex numbers up to sign in the following way: take a chart for one of the triangles it is on the boundary of, and take the difference of the two complex numbers that correspond to those two vertices in the chart. This number, up to sign, is the \emph{period} of the saddle connection. A collection of periods is a system of \emph{period coordinates} if they may vary in an open subset of $\cc^*/\{\pm 1\}^d$ and determine the periods of all other saddle connections in the triangulation. A system of period coordinates can be found for a Zariski open subset.\\

\noindent In \cite{lhc} we described how to cover $QD(\mathcal{M}_{g,n})$ by a locally finite collection of compact sets (with the respect to the Hausdorff topology on $QD(\mathcal{M}_{g,n})$ as a complex analytic space) which map to convex compact sets in $\cc^d$ with respect to period coordinate embeddings. The triangulations may degenerate on boundary points of compact sets. If we restrict to the unit area locus $QD^1(\mathcal{M}_{g,n})$, these sets form a locally finite closed cover of $QD^1(\mathcal{M}_{g,n})$. We may equip $\mathcal{QD}^1(\mathcal{M}_{g,n})$ with a metric $d_E$, which is the largest metric that is less than or equal to the standard Euclidean metric on $\cc^d$ on our chosen system of charts. The metric $d_E$ depends on choices but its local Lipschitz class does not.\\

\noindent The \emph{Teichm\"uller geodesic flow} $\{g_t | t \in \rr\}$ is the dynamical system on $QD(\mathcal{M}_{g,n})$ by acting on systems of charts (and on period coordinates) by $g_t(a + bi) \mapsto e^ta + e^{-t}bi$ for all $t \in \rr$, and complex numbers $a + bi, a,b \in \rr$. $X_1,X_2$ are said to be on the same \emph{strongly stable leaf} for the Teichm\"uller flow if there is a continuous path from $X_1$ to $X_2$ consisting of paths in convex compact period coordinate patches which the real parts of period coordinates remain constant. This property is invariant under $g_t$. Let $d^{ss}(X_1,X_2)$ denote the length of the shortest such path with respect to the metric $d_E$.\\

\noindent Our main theorem is the following:

\begin{theorem}Let $K$ be a compact subset of $QD^1(\mathcal{M}_{g,n})$ and let $\lambda$ denote Lebesgue measure on $\rr$. Let $\theta \in (0,1)$. There are positive real constants $C(K), r_0(K), \alpha(K,\theta)$ such that the following holds: if $T > 0,$ and $X, g_T(X) \in K$ and moreover, $$\lambda \{t \in [0,T]: g_t(X) \in K\} > \theta T, ~ \mathrm{and}$$ $$d^{ss}(X_1,X),d^{ss}(X_2,X) < r_0(K), ~ \mathrm{then}$$ $$\frac{d^{ss}(g_T X_1, g_T X_2)}{d^{ss}(X_1,X_2)} < C(K)e^{-\alpha(K,\theta)T}.$$ \end{theorem}

\noindent We will restate and prove this in a slightly different, but equivalent form, as \Cref{final}.\\

\noindent This is essentially a CAT(-1) type property. Of particular interest is how the constant $\alpha$ depends on $K$ and $\theta$. It is known that as $K$ exhausts $QD(\mathcal{M}_{g,n})$, or as $\theta \to 0^+$, the constant $\alpha$ becomes smaller. In particular, $g_T$ is not an Anosov flow. Previous results by \cite{abem}, using Hodge theoretic techniques, obtained a similar result, but required stronger hypotheses. In particular, they required that the associated Teichm\"uller flow segment spent a definite fraction of time on quadratic differentials with no short saddle connections; i.e. singularities of the flat metric were not allowed to come close to colliding. We summarize this as ``there is no loss of hyperbolicity near the multiple zero locus."\\

\noindent As applications, we prove two closing lemmas for the Teichm\"uller flow in the final section; one closing lemma for affine invariant manifolds, and the other for the full moduli space.\\

\subsection{Sketch of Proof}

\noindent One may associate a canonical choice of triangulation to almost every point along almost every Teichm\"uller flow line, which is dual to a pair of train tracks. There is a discrete set of times at which the triangulation changes, one dual train track splits along forward flow; the other splits along backward flow. In the context of periodic Teichm\"uller geodesic flow lines, this construction was discovered by Agol in \cite{Veering} and Gu\'eritaud in \cite{FG}. Their interest was in studying ideal traingulations of mapping tori associated to pseudo-Anosov homeomorphisms. Their choices of triangulation were different on the level of flat surfaces, but gave rise to the same triangulation of the mapping torus. We follow the construction in \cite{FG}.\\

\noindent The triangulation coarsely determines the half-translation surface under the additional assumption that the surface has no short essential simple closed curves. Having thus combinatorialized flow lines, we show that if a flow line spends enough time in a compact set, then pseudo-Anosov ``words" appear in its splitting sequence linearly often (see \Cref{ThicketySplit} and \Cref{pseudofilling}), and we prove a key property of pseudo-Anosov words composed sufficiently many times with a bounded amount of random ``noise" interspersed (as \Cref{CoroGoodWord}). The key property is that they induce contractions in Hilbert metrics on cones which correspond to subsets of $\mathcal{PMF},$ the space of measured foliations under the equivalence generated by scaling transverse measures and Whitehead moves.\\

\noindent A careful application of the pigeonhole principle shows that if a Teichm\"uller geodesic segment recurs enough to a compact set $K$, then such a sequence appears with a frequency proportional to the length of time we follow the geodesic. The key property is that given such a splitting sequence, after a bounded number of additional splits (which will not interfere with future splits along the way) the cone of equivalence classes of tangential measures on these train tracks is exponentially contracted in the appropriate Hilbert metrics. It is absolutely essential that we use tangential measures instead of the much smaller cone of transverse measures to deal with geodesics that spend most of their time near the multiple zero locus.\\

\noindent Finally, we compare the Hilbert metrics to the Euclidean metric $d^{ss}$.

\subsection{Acknowledgments}

\noindent The author would like to thank Alex Eskin, Simion Filip, Ilya Gekhtman, Ursula Hamenst\"adt, Sebastian Hensel, Howard Masur, Kasra Rafi, Saul Schleimer, and Robert Tang for helpful conversations.\\

\noindent Many of the constructions in section 3 were independently discovered by Vaibhav Gadre and Saul Schleimer, in an ongoing study of automata for producing pseudo-Anosov homeomorphisms compatible with given strata. Software by Mark Bell has been used to generate these automata for small strata of abelian and quadratic differentials. In ongoing work, Vincent Delecroix and Saul Schleimer are working to improve the speed of these computations and integrate them with other flat geometry software packages.\\

\noindent The author was supported by a Fields Postdoctoral Fellowship at the Fields Institute for Research in the Mathematical Sciences.\\

\includegraphics{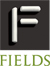}

\section{Teichm\"uller Space and Quadratic Differentials}

\noindent Let $g, n \geq 0$ be such that $3g-3+n > 0$. Let $\mathcal{M}_{g,n}$ be the moduli space of genus $g$ Riemann surfaces with $n$ distinct, unordered marked points. Let $\mathcal{T}_{g,n}$ be its orbifold universal cover, the Teichm\"uller space of genus $g$ Riemann surfaces with $n$ marked points. Then $\mathcal{T}_{g,n}$ is a complex manifold of dimension $3g-3+n$, homeomorphic to $\rr^{6g-6+2n}$. Let $S$ be a fixed closed genus $g$ surface with $n$ points deleted. Points in $\mathcal{T}_{g,n}$ consist of equivalence classes of the following: A compact Riemann surface with $n$ marked points $(X; x_1,...,x_n)$ and a homeomorphism $f: X \setminus \{x_1,...,x_n\} \to S.$ The equivalence relation is the following: $(X;x_1,...,x_n; f)$ and $(Y; y_1,...,y_n; g)$ are equivalent if there is a biholomorphism $h: X \setminus \{x_1,...,x_n\} \to Y \setminus \{y_1,...,y_n\}$ such that $g \circ h \circ f^{-1}$ is homotopic to the identity map on $S$. The homotopy class of the map $f$ is called a \emph{marking} of $X$.\\

\noindent This space can be given a metric as follows: If $K > 1$, then $(X;x_1,...,x_n;f)$ and $(Y;y_1,...,y_n;g)$ are distance at most $\log(K)/2$ if there is a $K$-quasiconformal map $h$ such that $g \circ h \circ f^{-1}$ is homotopic to the identity on $S$. This metric is the \emph{Teichm\"uller metric};  For precise characterizations of $K$-quasiconformal maps, and proof that this indeed forms a metric, see for instance \cite{Hubbard}.\\

\noindent Let $X$ be a genus $g$ Riemann surface and let $x_1,...,x_n$ be distinct points in $X$. Let $\mathcal{K}$ be a divisor whose associated bundle is the canonical bundle of $X$. The cotangent bundle of $\mathcal{T}_{g,n}$ at the point $(X; x_1,...,x_n; f)$ consists of the holomorphic sections of the degree $4g-4+n$ line bundle associated to the divisor $2\mathcal{K} + [x_1] + ... + [x_n]$; we refer to such sections as \emph{quadratic differentials}.\\

\noindent We may view them as meromorphic sections of the tensor square of the cotangent bundle, and the only poles allowed are simple poles at the points $x_i$. For a local coordinate $z$ on $X$ we have $q = f(z)dz^2$, and if $f(z) \neq 0,\infty$ we can find a local holomorphic coordinate $w$ s.t. $f(z)dz^{\otimes 2} = dw^{\otimes 2}.$ This holomorphic coordinate is unique up to rotation by $\pi$ and translations, and $\sqrt{q}$ is a $1$-form that can be integrated to give a local chart into $\cc$. Since transition maps between charts have derivative $\pm 1$, we can give $X$ a metric for which integrating $\sqrt{q}$ gives isometric charts. Such a metric will have a cone-type singularity of angle $(2 + k)\pi$ at $p$ if $q = f(z)dz^{\otimes 2}$ is such that $f$ vanishes to order $k \geq -1$ at $p$. We will refer to the zeros of $q$ and the points $x_1,...,x_n$ as \emph{singularities}, even if they happen to have cone angle $2\pi$. (This only happens at the points $x_i$, and for almost every $q$, they will have cone angle $\pi$.) If $q$ is a quadratic differential on $X$, it follows that $(X,q)$ admits a canonical path metric, for which geodesic arcs are Euclidean straight line segments than may turn only at cone points, and only such that the angle is at least $\pi$.\\

\noindent These transition maps also preserve the slopes of arcs, and a (possibly singular) area form $dx \wedge dy,$ where $x$ and $y$ are the real and imaginary parts of the local holomorphic coordinate function $w$ obtained by integrating a square root of $q$, i.e. a locally defined $1$-form $w$ such that $dw^{\otimes 2} = (dx + i dy)^{\otimes 2} = q$.

\noindent In summary, $(X,q)$ comes with the following structure: a flat metric with cone points, and for each slope (i.e. each element of $\rr\mathbb{P}^1$), a foliation whose tangent vectors at each point have that slope, that is allowed to have singularities at the cone points. (The point at infinity just corresponds to the vertical foliation.) The surface $X$, together with its collection of cone points and quadratic differential, is called a \emph{half-translation-surface}. Given a compact surface with a flat metric and a finite number of cone points whose holonomy from parallel transport about arcs not passing through cone points is $\pm I$, we can recover a Riemann surface structure and quadratic differential inducing that singular flat metric.\\

\noindent Classical deformation theory says that quadratic differentials are the cotangent space to Teichm\"uller space. The area of $(X,q)$ in the $q$-metric is the norm of $q$ for this cotangent bundle for the Teichm\"uller metric on Teichm\"uller space; the dual norm on the tangent bundle is a Finsler norm for which the 1-parameter families $\{ g_t(X,q): t \in \rr \}$ are geodesics and the parameter $t$ is distance.\\

\noindent We generally will not consider the everywhere zero quadratic differential on any Riemann surface, and so the complement of the zero section of the cotangent bundle of $\mathcal{T}_{g,n}$ will be referred to as the Teichm\"uller space of quadratic differentials, and we shall denote it by $QD(\mathcal{T}_{g,n})$. If we want to quotient by the action of the mapping class group, we write $QD(\mathcal{M}_{g,n}).$ We will write $QD^1$ instead of $QD$ if we want to restrict to half-translation surfaces of unit area. By the uniformization theorem, there is a hyperbolic (complete, sectional curvature everywhere equal to $-1$) defined on the complement of $n$ points in a compact Riemann surface of genus $g$ (provided $g > 0$.) If $X = (Y;y_1,...,y_n;f)$ defines a point in Teichm\"uller space, then we will write $\sigma_X$ to denote this metric on $Y \setminus \{y_1,...,y_n\}$). Sometimes we will abuse notation and use the same letter to denote the Riemann surface and a point in Teichm\"uller space.\\

\noindent A \emph{saddle connection} is the image of a locally constant speed geodesic map $\gamma: [0,1] \to X$ with the property that $\gamma(0)$ and $\gamma(1)$ are the only singularities in the image of $\gamma$. We do not require $\gamma(0)$ and $\gamma(1)$ to be distinct singularities. $X$ can be given a CW structure whose 0-skeleton consists of cone points, whose 1-skeleton consists of saddle connections, and whose open 2-cells are isometric to triangles in Euclidean space. We will call any such a CW structure a \emph{triangulation} of the surface, even though they are not necessarily simplicial complexes: \emph{we require each triangle in a triangulation to have three distinct edges, but do not require that it have three distinct vertices. We also require that every singularity be a 0-cell of the triangulation.} We describe one way of finding a triangulation in the next section.\\

\noindent Fix a triangulation of $(X,q)$ by saddle connections. Assume $q$ has is in the \emph{principal stratum}, that is to say, $q$ has $4g-4 + n$ simple zeros and $n$ simple poles. Then there is a collection of $6g - 6 + 2n$ saddle connections in the triangulation such that the complex lengths of the edges $\gamma_j, 1 \leq j \leq 6g - 6 + 2n$, defined up to sign by $\int_{\gamma_j} \sqrt{q}$ locally determine $(X,q)$. These complex lengths are called \emph{periods} of the saddle connections, and for each point in the principal stratum, some collection of periods forms a local holomorphic coordinate system for $QD(\mathcal{T}_{g,n})$. Such coordinate systems are called \emph{period coordinates.} We will refer to the absolute values of the real and imaginary parts of saddle connections as \emph{widths} and \emph{heights}. The complement of the principal stratum is an analytic subvariety, so it is closed and nowhere dense.\\

\noindent The $6g-6+2n$ dimensions can be interpreted as a cohomology group. There is a degree 2 branched cover $\tilde{X}$ of $X$ branched at the zeros and poles of $q$ such that the pullback of $q$ to $\tilde{X}$ is the square of a holomorphic $1$-form; the cohomology class of this $1$-form necessarily lies in the subspace $H_{odd}^1(\tilde{X};\cc) \subset H^1(\tilde{X};\cc)$ consisting of $-1$-eigenvectors of the action of the involution $\iota$ of $\tilde{X}$ interchanging the two sheets of the branched cover. A basic dimension count, using Riemann-Hurwitz formula for $\chi(\tilde{X})$, shows that the $1$-eigenvectors of $\iota$ are exactly the image of $H^1(X)$ under pullback. By naturality of the cup product, these two eigenspaces are symplectically orthogonal to each other.\\

\section{The $L^\infty$ Delaunay Triangulation and Associated Train Tracks}

\begin{definition} A \emph{train track} is a finite ribbon graph $\tau$ embedded in a surface with a well-defined tangent line at each vertex, and such that the angles formed by consecutive outward-pointing half-edges from each vertex are either $0$ or $\pi$. Moreover, at each edge there is at least one angle of $0$ and there are exactly two angles of $\pi$.\end{definition}

\begin{definition}A \emph{train route} is a sequence of oriented edges $e_1,...,e_k$ of a train track $\tau$, possibly with repetition and such that an edge need not have the same orientation every time it appears, with the property that there is an immersion $\gamma$ of $[0,k]$ into the surface with the property that $\gamma$ maps $[(r-1),r]$ onto $e_r$ with the designated orientation for $1 \leq r \leq k$. \end{definition}

\noindent We will usually assume our train tracks are trivalent, or \emph{generic}. At each vertex $v$ of a generic track $\sigma$, for each sufficiently small $\epsilon > 0$ there are three smooth arcs of length $\epsilon$ in $\sigma$ starting at $v$: two that are contained in \emph{small} half-edges, which leave $v$ in the same direction (same derivative up to positive scalar), and one \emph{large} half-edge, which leaves $v$ in the opposite direction from the other two. Vertices of a train track are called \emph{switches}. An edge that is large at both ends is said to be \emph{large}; an edge that is small at both ends is \emph{small}, and an edge that is large at one end and small at the other is \emph{mixed}.\\

\begin{definition} A \emph{transverse measure} is a collection of nonnegative real numbers, one for each edge and not all zero, satisfying the \emph{switch condition}: at each vertex $v$, the sum of the weights of the edges of $v$ in one direction is equal to the sum of the weights on the edges leaving $v$ in the opposite direction.\end{definition}

\noindent For a generic track, this means the large edge has weight equal to the sum of the sum of the weights of the two small edges at each vertex.\\

\noindent The complement of a train track $\tau$ in a surface is a union of open subsurfaces. Let $D$ be such a subsurface, and assume $D$ is homeomorphic a disk or a disk with punctures.

\begin{definition}A \emph{side} of the complementary region $D$ is a maximal train route such that for each pair of consecutive oriented edges $e_i$ and $e_{i+1}$ meeting at vertex $v_i$ in the train route, $e_{i+1}$ is the successor of $v$ in the circular counterclockwise order of the edges about $v_i$.\end{definition}

\noindent For $n=1,2,3$ we use the terms monogon, bigon, trigon. By this we mean that there are $n$ points at which we cannot choose the inward normal vector continuously; call a maximal arc on $\partial D$ for which the inward normal can be chosen continuously a \emph{side} of $D$; call the points of discontinuity in the normal the \emph{vertices} of $D$. At each vertex of $D$, the two sides that meet (possibly the two ends of the same side if $D$ has only one side) form an internal angle of $0$. When a train track is embedded on in a half-translation surface, we will mainly consider train tracks with the following property: each complementary region is a disk containing exactly one cone point, and it is an $n$-gon if the cone point has cone angle $n \pi$, and no cone point lies on the train track.

\begin{definition}
Assume every complementary region of $\sigma$ is an $n$-gon, where $n \geq 1$ is allowed to depend on the region. A \emph{tangential measure} on $\sigma$ is a collection of positive weights with the property that for each $n$-gon $D$ without punctures, the sum of the weights on the boundary of $D$ is more than twice the weight of each side of $D$. Two tangential measures are \emph{equivalent} if their difference is in the span of vectors of the form $E_v - F_v$, where $E$ is assigns weight 1 to all edges on one side of a vertex $v$ and $F$ assigns weight $1$ to all edges on the other side of $v$.
\end{definition}

\noindent For example, if $a$ is large at a vertex $v$ and $b$ and $c$ are small at $v$, we can add $t$ to the weight of $a$ and $-t$ to the weights of $b$ and $c$ and remain in the same equivalence class, provided the original tangential measure assigned weight at least $t$ to each of $b$ and $c$.

\begin{definition} A train track $\tau$ on $S_g \setminus \{x_1,...,x_n\}$ is \emph{filling}, or \emph{fills} $S_g \setminus \{x_1,...,x_n\},$ if every complementary region of $\tau$ is a disk or once-punctured disk.\end{definition}

\begin{definition} The train track $\tau$ is \emph{complete} if every complementary region is a trigon or once-punctured monogon.\end{definition}

\begin{definition} A filling train track $\tau$ on a surface of genus $g$ with $n$ marked points is \emph{birecurrent} if it satisfies the following:\begin{itemize} \item No marked point is on $\tau$ \item Every complementary region $R$ of $\tau$ is an $m_R$-gon for some $m_R \geq 1$ \item Every complementary region not containing a marked point is an $m$-gon for some $m \geq 3$ \item $\tau$ admits a transverse measure with strictly positive weight on every edge \item $\tau$ admits a tangential measure with strictly positive weight on every edge.
\end{itemize} \end{definition}

\begin{definition}An \emph{oriented rectangle} $R$ in a half-translation surface $(X,q)$ is the closure of a contractible open set $U \subset X$ with no singularities such that for any $z_0 \in U$ the developing maps $z \mapsto \int_{z_0}^z \sqrt{q}$ map $U$ homeomorphically onto an open rectangle in $\cc$ whose boundary has vertical and horizontal sides. We call $R$ an \emph{oriented square} if the vertical and horizontal sides have equal length. We say that $U$ is the \emph{interior} of $R$, and $R \setminus U$ is the boundary of $R$.\end{definition}

\noindent Remark: it might be the case that the interior and boundary of an oriented rectangle $R$ with interior $U$ are not the interior and boundary of $R$ in the usual point-set topological sense due to identifications of segments, but they will be the interior and boundary of $U$.\\

\noindent We can build Riemann surfaces with quadratic differentials from filling birecurrent train tracks with transverse measure $\mu$ and tangential measure $\nu$ as follows. Suppose that $\tau$ is a generic train track on a closed oriented surface $S$, equipped with transverse and tangential measures that assign positive weight to each edge. Moreover, assume each complementary region $R$ of $\tau$ is an $m_R$-gon, where each $m_R$ is \emph{odd}. (We do not require $m_{R_1} = m_{R_2}$ for distinct regions $R_1$ and $R_2$.) Then, there is a unique half-translation structure on $S$ which is a union of oriented rectangles with disjoint interiors. The widths of the rectangles are the weights of the transverse measure and heights are the weights of the tangential measure.\\

\noindent The switch conditions and the ribbon structure tell us how to isometrically identify the horizontal edges on each side of each switch - the side of the large half-edge is glued to the other two, and the clockwise order a small circle visits the three rectangles is the same as the clockwise order at which a small circle about the corresponding vertex visits the three edges. Typically there is no singularity at the point these three rectangles meet.\\

\noindent For each complementary $3$-gon $D$ with side $s$, the aforementioned identifications of horizontal edges concatenate vertical edges of the rectangles coming from an edges of $\sigma$ on $s$. If these three concatenated vertical segments have total weights $a+b,b+c,a+c,$ where $a,b,c > 0$ we can identify them and attain a cone point of cone angle $3\pi$ that is distance $c,a,b$ away from their shared vertices. For each complementary $1$-gon we take each side made of concatenated edges and identify points equidistant from its midpoint, where there is a cone point of angle $\pi$.

\begin{definition}\label{RectDecomp}We refer to the above construction as the \emph{rectangle decomposition} associated to $\mu$ and $\nu$. \end{definition}

\noindent The same construction is valid if we allow $n$-gons where $n \neq 1,3$ but there is additional choice involved for each $n$-gon with $n$ even.\\

\noindent Notice that replacing the tangential measure $\nu$ by an equivalent one does not actually change the underlying Riemann surface or meromorphic quadratic differential. Since the dimension of $QD(\mathcal{T}_{g,n})$ is $6g-6+{2n}$, basic Euler characteristic and dimension counting implies that any small perturbation of transverse or tangential measure that is not just a replacement of $\nu$ by an equivalent measure changes the quadratic differential.\\

\begin{definition}We say that $(X,q)$ is \emph{suited} if its vertical and horizontal foliations do not have any saddle connections as leaves.\end{definition}

\noindent It is well known that if the foliation associated to a given slope does not contain any saddle connections, then every leaf is dense.\\

\noindent We now describe how to choose a ``best" triangulation and an associated pair of train tracks, given a generic quadratic differential, i.e. a suited quadratic differential such that the real part of the period of every saddle connection is different from the imaginary part of the period of every saddle connection. (This fails on a set that is dense but Masur-Veech measure zero.) This condition can be relaxed somewhat; we only need to consider saddle connections that are diagonals of oriented rectangles and which are not too long compared to the diameter of the surface with the $q$-metric, but anyway, we can stick to our stronger assumptions. Later we will see (\Cref{greedy}) that it is sufficient to check finitely many inequalities to tell if a triangulation is the best, and the choice of a pair of best train tracks therefore makes sense on a dense open set.\\

\noindent \begin{definition} Define the flat $L^\infty$ metric on $(X,q)$ to be the metric for which a $q$-geodesic line segment $\gamma$ with $$\int_\gamma \sqrt{q} = \pm(a+bi)$$ \noindent has length $\max(|a|,|b|)$, or more generally, such that the length of any curve $\gamma$ is $$\max\left(\int_\gamma |\mathrm{re}(\sqrt{q})|, \int_\gamma |\mathrm{im}(\sqrt{q})|\right).$$\end{definition}

\noindent We note that if $(X,q)$ has no marked points, then its universal cover $(\tilde{X},\tilde{q})$ is nonpositively curved in the sense of Alexandrov. Unfortunately, if $(X,q)$ has marked points, this is no longer the case, since there is positive curvature at the poles of $q$. We will instead build a complete nonpositively curved branched cover of $(X,q)$, and assume that is branched only over the marked points, with local degree at least $2$ at each ramification point. In particular, this cover will have the property that there is a unique geodesic joining each pair of points.\\

\noindent If $n$ is even, we can fix a collection of smooth arcs whose endpoints are a partition of the marked points into pairs, and declare these arcs to be branch cuts for a 2-sheeted branched cover. If $n$ is odd, we can take any double cover of the compact surface $X$ and then apply the procedure for $n$ even. Only finitely many covers can be obtained in this way, since a finitely generated group has only finitely many index 2 subgroups. In any case, we can go a step further and pass to a normal (Galois) cover of degree at most 24. We can then pull back the quadratic differential to the universal cover of this finite cover of the compact surface; let $(\tilde{X},\tilde{q})$ be this infinite branched cover of $(X,q)$.\\

\noindent The leaves of the vertical and horizontal foliations are convex in the $\tilde{q}$-metric - even the leaves containing cone points are convex, and the transverse measures on the foliations induce metrics on the leaf space. (The leaf space of a measured foliation is known to be an $\rr$-tree, i.e. a space where there is a unique path between each pair of points, realized by the projection of any $\tilde{q}$ geodesic with endpoints on the given leaves.) The $L^\infty$ flat distance between two points $x$ and $y$ is just the max of the distance between their vertical leaves and the distance between their horizontal leaves, since $\tilde{q}$-geodesics project to unparametrized geodesics in the leaf space of $(\tilde{X},\tilde{q})$. That is to say, if $a,b,c$ are three points on a $\tilde{q}$-geodesic with $b$ between $a$ and $c$, and $d_h$ is the transverse distance on the space of vertical leaves, then $$d_h(a,b) + d_h(b,c) = d_h(a,c).$$ The same is true for transverse distance on the space of horizontal leaves.

\noindent \begin{definition} On $(\tilde{X},\tilde{q})$ of $(\hat{X}, \hat{q}),$ We declare the \emph{singularities} of $(\tilde{X},\tilde{q})$ to be any points that project to singularities of $(X,q).$ We will refer to $(\tilde{X},\tilde{q})$ as a \emph{normal pole-free-disk cover} of $(X,q)$; it is a complete, simply connected Alexandrov non-positively curved surface, and the collection of isometries that that commute with the projection $\pi: (\tilde{X},\tilde{q}) \to (X,q)$ is transitive on $\pi^{-1}(p)$ for each $p \in X$; we will refer to this as the \emph{deck group} of the normal pole-free-disk cover.\end{definition}

\noindent A normal pole-free-disk cover of $(X,q)$ can be given a $CW$ structure whose open 2-cells correspond to the cone points: for each point $p \in \tilde{X}$, if there is a unique closest cone point $z$ in the $L^\infty$ flat metric, $p$ belongs to the open 2-cell corresponding to $z$. This cellulation of $\tilde{X}$ is the \emph{Voronoi diagram} and it is invariant under the deck group of the normal pole-free-disk cover.\\

\noindent No point $p$ can belong to the boundary of more than 3 cells, since if it did, then the maximal square disk centered at $p$ with no cone points in its interior would have four cone points on its boundary, and there would be a horizontal or vertical saddle connection, or a pair of saddle connections, one of which has width equal to the other's height, which we have disallowed. If $p$ belongs to the boundary of exactly three cells, draw a triangle of saddle connections joining them. These triangles project down to the triangles of the $L^\infty$ Delaunay triangulation of $(X,q)$. This triangulation was previously considered in \cite{FG}, where it was referred to as the Delaunay cellulation. This triangulation has the following, alternative characterization: find all embedded, singularity-free open squares, which are maximal as oriented squares (i.e. oriented rectangles not properly contained in any other oriented squares). Connecting them by saddle connections produces all of the edges, and each triangle will have all of its vertices on one of the squares. In particular, each triangle has at least one edge with positive slope and at least one edge with negative slope.\\

\noindent In \cite{FG} it is assumed that there are no vertical or horizontal saddle connections in the $L^\infty$ Delaunay triangulation, though a square with one vertex on each side is permitted, in which case the cellulation thus obtained may have quadrilaterals. However, this can be resolved by applying forward or backward Teichm\"uller flow by a sufficiently small time; see \Cref{Successor}.\\

\noindent Our choice of name for this triangulation is in part to distinguish it from the Delaunay triangulation used in \cite{Delaunay}, in which essentially the same construction was done with the usual ($L^2$) flat metric associated to $q$ instead of the $L^\infty$ flat metric. The Delaunay triangulation can be characterized by the property that in any normal pole-free-disk cover, the circumcircles of triangles bound cone-point-free Euclidean disks.\\

\begin{proposition}\label{PropVertSeparatrix} Let $v$ be a cone point and $\gamma$ a vertical geodesic ray emanating from $v$. Then there is some maximal embedded square with three cone points on its boundary such that the inital segment of $\gamma$ is in the interior of the square, and such that $\gamma$ is between the other two boundary cone points. Conversely, for each maximal embedded oriented square with exactly three cone points on its boundary, one vertex is horizontally between the other two, and therefore on the top or bottom edge of the square. Thus there is a bijection between triangles and maximal embedded oriented squares with three boundary singularities. \end{proposition}

\noindent The ray $\gamma$ is called a \emph{separatrix}, or \emph{vertical separatrix}. We will refer to maximal embedded oriented squares with three singularities on the boundary as \emph{circumsquares} of $L^\infty$ Delaunay triangles.\\

\noindent Proof: To prove the first claim, we note that the diameter of the original surface is finite, so there do not exist arbitrarily large embedded cone-point-free squares. Thus, the squares bisected by the initial segment of $\gamma$ have bounded side length, and at the least upped bound another cone point appears on the boundary. If two cone points occur simultaneously we are done. If only one cone point $w$ occurs on the left or right side of the square, we consider the family of squares with $v$ on the top or bottom $w$ the left or right; eventually we will hit a new cone point, by bounded size of the squares, and in this case we are done. The final case is that $w$ occurs on the opposite side of the square from $v$; in this case, by symmetry we assume $v$ is on the bottom, $w$ is on the left half of the top. Then we may slide the square to the right while keeping $w$ and $v$ on the top and bottom of its boundary until $v$ reaches the corner or another cone point appears on the boundary. If a cone point appears we are done. If not, then $w$ is in the top-left corner. We now consider the family of squares with $v$ on the bottom and $w$ on the left side. As before, there is a first time at which a new cone point appears on the boundary. This completes the proof of the first claim. The rest of the proposition is immediate. $\Box$\\

\noindent We now describe a way to get a train track, and a choice of transverse and tangential measures, from a class of triangulations that include the $L^\infty$ Delaunay triangulation. First we introduce the class of triangulations.

\begin{definition}A triangulation $T$ by saddle connections is a \emph{veering triangulation} if in any normal pole-free-disk cover, every saddle connection of $T$ is the diagonal of an oriented rectangle.\end{definition}

\begin{proposition}\label{EquivSlope}Suppose $T$ is a triangulation of a normal pole-free-disk cover $(\tilde{X},\tilde{q})$ of $(X,q)$ by saddle connections, none of which is vertical or horizontal, with every cone point as a vertex, and $T$ is invariant under the deck group of the cover. Then $T$ is veering if and only if no triangle has three edges of positive slope or three edges of negative slope.\end{proposition}

\noindent Proof: If there is a triangle in the plane whose sides all have positive slope, then it is obtuse and the longest side is the diagonal of an oriented rectangle containing the third vertex. This is impossible for any triangle in a veering triangulation. The same is true for a triangle with all negatively sloped sides.\\

\noindent Now assume that there are cone points $a,b$ joined by a saddle connection $ab$ of $T$ that is not the diagonal of an oriented rectangle; then there is some cone point $c$ such that the vertical distances $d_v(a,c)$ and $d_v(b,c)$ are both less than $d_v(a,b)$, and the horizontal distances $d_h(a,c)$ and $d_h(b,c)$ are both less than $d_h(a,b)$. Over all possible choices of such $ab$ and $c$, we choose one such that the $q$-distance from $c$ to the segment $ab$ is minimal. Then there is a trapezoid with one vertical side, one horizontal side, $ab$ and a side parallel to $ab$ passing through $c$ with no cone points in its interior. No saddle connection of $T$ has initial segment inside the acute angles $bac$ or $abc$, because such a saddle connection would pass closer to $c$ than $ab$ does. If $ab$ is extended to a bi-infinite geodesic $\ell$, then there is a triangle of $T$ with vertices $a,b$ and a third vertex $d$ in the connected component of $\tilde{X} \setminus \ell$ that contains $c$. Suppose $c \neq d$. If $ad$ and $bd$ are extended to bi-infinite geodesics not containing $c$, then $c$ is on the same side of each of the infinite geodesics $ab$, $ac$, and $ad$ as the interior of the triangle $abd$, so $c$ is in the interior of triangle $abd$, a contradiction. Thus $c = d$, and $abc$ is a triangle with 3 sides of the same slope. $\Box$

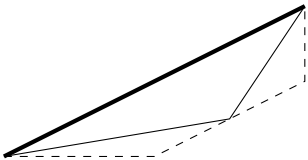
\begin{figure}[H]
\begin{tikzpicture}
\draw [dashed] (0,0) --(2,0) -- (4,1) -- (4,2);
\draw (0,0) --(3,0.5) --(4,2);
\draw [ultra thick] (0,0) --(4,2);
\end{tikzpicture}
\caption{The proof of \Cref{EquivSlope}. The thick edge of triangle $abc$ is $ab$. The remaining three sides of the trapezoid are dotted, and $c$ is on the side of the trapezoid parallel to $ab$.}
\end{figure}

\begin{proposition}\label{PropCanonicalRectangles} Assume $(X,q) \in QD(\mathcal{T}_{g,n})$ is suited, and equipped with a veering triangulation $T$ by saddle connections. For each triangle in a normal pole-free-disk cover $(\tilde{X},\tilde{q})$, let $v$ be the vertex that is horizontally in between the other two vertices, $a$ and $b$. Let $\gamma^!$ be the geodesic segment that is horizontal, passes through the midpoint of the saddle connection joining $a$ and $b$, and whose endpoints are on the vertical leaves of $a$ and $b$. Let $\gamma$ be the vertical segment from $v$ to the horizontal $\gamma^!$.\\

\noindent Then, the complement of the segments $\gamma$ and $\gamma^!$ generated by the above procedure is a disjoint union of interiors of oriented rectangles, and it is invariant under the group of covering transformations for the branched cover $\tilde{X} \to X$.\end{proposition}

\noindent Proof: The second claim is obvious, given the first. We will refer to the vertical segments as \emph{drawn separatrices}, and the horizontal segments as \emph{separatrix stoppers}. We break the proof into four steps.\\

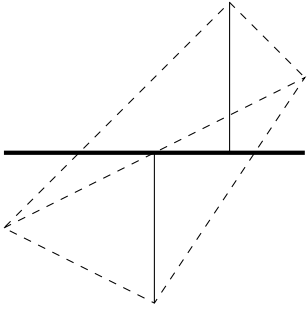
\begin{figure}[H]
\begin{tikzpicture}
\draw [dashed] (0,1) --(4,3) -- (2,0) -- (0,1) -- (3,4) -- (4,3);
\draw (2,0) --(2,2);
\draw (3,4)--(3,2);
\draw [ultra thick] (0,2) --(4,2);
\end{tikzpicture}
\caption{Two triangles of the triangulation (dashed) and the corresponding drawn separatrices (solid thin segments) and their shared separatrix stopper (solid thick segment). The shared edge of the triangles is dual to a large edge of a train track, since it is the widest edge of both triangles.}
\end{figure}

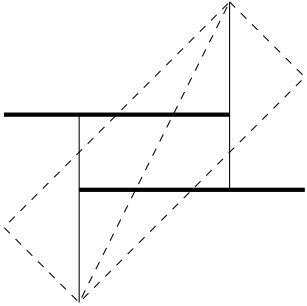
\begin{figure}[H]
\begin{tikzpicture}
\draw [dashed] (3,4) --(4,3) -- (1,0) -- (0,1) -- (3,4) -- (1,0);
\draw (3,4) --(3,1.5);
\draw (1,0) --(1,2.5);
\draw [ultra thick] (1,1.5) --(4,1.5);
\draw [ultra thick] (0,2.5) --(3,2.5);
\end{tikzpicture}
\caption{Two triangles of the triangulation (dashed) and the corresponding drawn separatrices (solid thin segments) and a separatrix stoppers (solid thick segments). The shared edge of the two triangles is dual to a small edge of the train track, since it is not the widest edge of either triangle.}
\end{figure}

\begin{figure}[H]
\begin{tikzpicture}
\draw [dashed] (0,2) --(1,0) -- (3,1) -- (0,2) -- (4,4) --(3,1);
\draw (1,0) --(1,1.5);
\draw (3,1) --(3,3);
\draw [ultra thick] (0,1.5) --(3,1.5);
\draw [ultra thick] (0,3) --(4,3);
\end{tikzpicture}
\caption{Two triangles of the triangulation (dashed) and the corresponding drawn separatrices (solid thin segments) and a separatrix stoppers (solid thick segments). The shared edge of the triangles is dual to a mixed edge of the train track, since it is the widest edge of only one triangle.}
\end{figure}
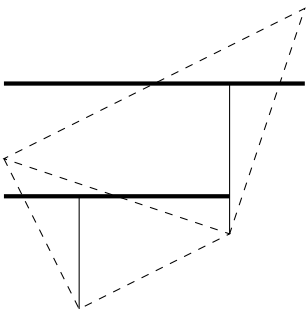

\noindent Throughout the proof, we will often refer to a rectangle $R$ in $\tilde{X}$ and say that points and sets are to the left or to the right of each other. What this means that we will fix a choice of up down, left, and right in $R$ that is compatible with $\tilde{q}$, and we can choose a (bi-infinite) vertical or horizontal geodesic $g$ in $\tilde{X}$ passing through $R$ or $\partial R$ that divides $\tilde{X}$ into three disjoint pieces: $g$, and two connected components $A$ and $B$ of $X - g$. For some choice of $g$, $A$, $B$, we may say $P$ is above/below/to the left of/to the right of $Q$ if $P$ is contained in $g \cup A$ and $Q$ is contained in $B$.\\

\noindent STEP 1: Every end of a drawn separatrix is a separatrix stopper.\\

\noindent Proof: This is clear.\\

\noindent STEP 2: Each end of each separatrix stopper is a non-end point of a drawn separatrix.\\

\noindent Proof: Let $\gamma,\gamma^!,v,a,b$ be as in the proposition. The union of the cone-point-free rectangles associated to the three edges of the triangle they form is a larger cone-point-free rectangle $R$. What we need to show is that the drawn separatrices emanating from $a$ and $b$ are long enough to meet the ends of $\gamma^!$. There are cases depending on which of $a,b,v$ lie on which sides of $R$, but since we may rotate by $\pi$ and reflect across vertical or horizontal lines, we may assume $v$ is on the bottom, $a$ is on the left, and $b$ is at the top-right corner.

\noindent First, we need to show that the upward drawn separatrix $\hat{\gamma}$ from $a$ on the boundary of our square has length at least half of the height of the saddle connection $\overrightarrow{ab}$ from $a$ to $b$. If it is not, then the triangle associated to this separatrix contains two saddle connections starting at $a$ and of smaller height than the saddle connection $\overrightarrow{ab}$. Pick $\epsilon > 0$ the be shorter than the altitude of any triangle in $T$, and draw a circle $C$ of radius $\epsilon$ about $a$. Starting from the first point on $\hat{\gamma}$ that meets $C$, travel along $C$ counterclockwise. The first saddle connection of $T$ that $C$ crosses will be one of the two sides of the triangle in $T$ that contains the initial segment of $\hat{\gamma}$; another such side is obtained by moving clockwise along $C$ from the same starting point. The cone-point-free rectangle whose diagonal is $\overrightarrow{ab}$ guarantees that any such saddle connection has height at least as large as $\overrightarrow{ab}$, so the average height of the two saddle connections is more than half the height of $\overrightarrow{ab}$. This implies that the drawn segment of $\hat{\gamma}$ passes through the left endpoint of $\gamma^!$. Similarly, the rectangle with diagonal $\overrightarrow{bv}$ guarantees that the downward drawn separatrix from $b$ must have height at least half as large as the height of $\overrightarrow{bv}$, which is greater than half the height of $\overrightarrow{ba}$.\\

\noindent STEP 3: If a drawn separatrix intersects a separatrix stopper, the intersection point is at the end of either the separatrix or the separatrix stopper, but not both.\\

\noindent Proof: If both end at the same point then $(X,q)$ has a vertical or horizontal saddle connection, which we don't allow. Otherwise, consider the saddle connection $s$ whose midpoint is the midpoint of the separatrix stopper. Otherwise, let $a,b,v, \gamma, \gamma^!$ be as before: $R$ is a cone-point-free open rectangle with $v$ on the bottom and $a$ on the left, and $b$ in the top-right corner. Let $v^\prime$ be the cone point from which a drawn separatrix $\gamma^\prime$ emanates, and assume $\gamma^\prime$ crosses $\gamma^!$. Let $a^\prime,b^\prime$ be the other two vertices of the triangle that contains the initial segment of $\gamma^\prime$, and let $R^\prime$ be the oriented rectangle with $a^\prime, b^\prime, v^\prime$ all on its boundary. We have two cases, depending on whether $\gamma^\prime$ enters $R$ from above or from below.\\

\noindent Case 1: $\gamma^\prime$ enters $R$ from below. Then the rectangle $R^\prime$ associated to $\gamma^\prime$ cannot contain $v$, so $R^\prime \cap R$ lies entirely to the left of $v$ or entirely to the right of $v$. Now, $a^\prime$ and $b^\prime$ cannot both be to the left of $R$ or both to the right of $R$, since the segments joining them to $v^\prime$ have opposite slopes. They also cannot both be below $R$, so it follows that at least one vertex of $R^\prime$ is above $R$, assume this vertex is $a^\prime$. Since $R^\prime$ cannot contain $a$ or $b$, $a^\prime$ cannot be to the left or to the right of $R$. But then $\overrightarrow{v^\prime a^\prime}$ must cross $\overrightarrow{ab}$, which is impossible.\\

\noindent Case 2: $\gamma^\prime$ enters $R$ from above. Since $\gamma^\prime$ crosses $\gamma^!$. Now, if $\gamma^\prime$ has a vertex below $a$, then that vertex must be below $R$, since it cannot be in $R$, and if it is to the left or right of $R$ then $R^\prime$ would contain $a$ or $b$. But then  $\overrightarrow{v^\prime a^\prime}$ or $\overrightarrow{v^\prime b^\prime}$ would cross $\overrightarrow{ab}$, which is impossible. So both edges are above $a$ or at the same height, in which case one edge is $a$ because there are no horizontal saddle connections. Since the midpoint of $\overrightarrow{a^\prime b^\prime}$ is lower than the midpoint of $\overrightarrow{ab}$ it follows that both $a^\prime$ and $b^\prime$ are below $b.$ Since $R^\prime$ cannot contain $b$, neither $a^\prime$ nor $b^\prime$ can be to the right of $b$. Therefore the triangle containing the initial segment of $\gamma^\prime$ contains a cone point which is to the right of $\gamma^\prime$, to the left of $b$, and below $b$, and above $a$. There is no such cone point, because $R$ is cone-point-free.\\

\noindent STEP 4: The full proposition.\\

\noindent Proof: Every point $p \in \tilde{X}$ belongs to some triangle of $T$, so if $p$ is not on a drawn separatrix or separatrix stopper there is a vertical segment from $p$ that hits this separatrix stopper, and the length of the shortest such segment is bounded by the side length of the largest square. Now, given a separatrix stopper $\gamma^!$, a point $x$ in $\gamma^!$ that is not part of a separatrix, and a continuous choice $N$ of normal direction to $\gamma^!$, we consider the maximal connected subset $S$ of the separatrix stopper that contains $x$ with the property that for all sufficiently small $\epsilon >0$, the set of vertical segments of length $\epsilon$ leaving $S$ in direction $N$ with all parts of drawn separatrices deleted is connected. $S$ is an interval; let $S^\circ$ be $S$ with its endpoints deleted.\\

\noindent We claim that $S$ is the edge of a rectangle. There must be some distance after which a trajectories leaving $S^\circ$ in direction $N$ intersects a drawn separatrix or separatrix stopper, by boundedness of lengths of vertical segments that avoid separatrix stoppers. Let $d$ be the infimum such distance. $d$ be positive and must actually be attained, since there are locally finitely many drawn separatrices and separatrix stoppers. It is impossible for vertical trajectory to hit a separatrix without hitting a separatrix stopper unless it is already on the separatrix.\\

\noindent Now, suppose that there is not a separatrix stopper extending all the way across the segment distance $d$ from $S^\circ$. Then, at either the left or right end of this separatrix stopper we contradict step 3, so this is not the case.\\

\noindent Thus the trajectories leaving $S^\circ$ form an open rectangle with an equal length segment at the other end. We must check that it is a connected component of the complement. By construction, the left and right vertical trajectories on the boundary of this rectangle begin as drawn separatrices. If one of these drawn separatrices ends at a distance less than $d$, it cannot do so by meeting a drawn separatrix, so the boundary of this rectangle passes through a cone point, and continues along another drawn separatrix, and this separatrix must continue until it reaches a separatrix stopper. This establishes that the union of the trajectories leaving $S^\circ$ in direction $N$ that continue until reaching separatrix stoppers is an open rectangle.$\Box$

\begin{proposition}
No two separatrix stoppers intersect.
\end{proposition}

\noindent Proof: First we rule out the case in which they intersect just at a shared endpoint. Say the endpoint lies on a separatrix emanating from $v$ and the other two endpoints are on separatrices emanating from $v^\prime$ and $v^{\prime\prime}$. Then $vv^\prime$ and $vv^{\prime\prime}$ are saddle connections in the triangulation, and their associated oriented rectangles share a boundary edge which contains part of a separatrix emanating from $v$, which contains no singularities. The heights of the saddle connections these separatrix stoppers bisect are equal, so the union of these two rectangles is a rectangle with no singularities on the interior, with $v$ on one horizontal edge, with with $v^\prime$ and $v^{\prime\prime}$ vertices of the other horizontal edge. This contradicts the absence of saddle connections.\\

\noindent Now, we rule out the case in which two separatrix stoppers intersect on a segment of positive length. We note that the saddle connections bisected by these separatrix stoppers cannot have equal height, since otherwise there would be a horizontal saddle connection joining them as in the previous case. So one saddle connection has strictly larger height. Since each saddle connection is the diagonal of an oriented rectangle, neither rectangle can contain an edge of the other. This means that the rectangle of smaller height must intersect both vertical edges of the rectangle with larger height. But then the diagonals of these two rectangles are a pair of saddle connections of the triangulation that cross each other, which is impossible.$\Box$

\begin{proposition}The rectangle decomposition on $\tilde{X}$ descends to $X$.\end{proposition}

\noindent Proof: The branched covering $\tilde{X} \to X$ is a covering map in the complement of the cone points, so the group of covering transformations acts on the union of the open rectangles by isometries without fixed points. Therefore each rectangle of $\tilde{X}$ injects to $X$. $\Box$\\

\noindent Now that we have established that there is a rectangle decomposition of $(X,q)$, we claim that there is an associated train track for this rectangle decomposition. Some pairs of separatrices share a separatrix stopper. We now perform the following operation to modify our rectangle decomposition: We draw horizontal lines a small distance $\epsilon > 0$ above and below this separatrix stopper, and erase what is drawn between them. The horizontal lines, together with the two drawn separatrices that pass through ends of the separatrix stopper form a rectangle. We thus make the four rectangles whose boundary intersected our separatrix stopper shorter by $\epsilon$, and create a new rectangle of height $2\epsilon$. Now, our system of rectangles with their heights and widths gives us a generic train track with a positive transverse measure and a tangential measure that depended on $\epsilon$ but whose equivalence class was well defined. The edges of the train track are rectangles, and the vertices are the points where separatrices end; the large edge at each vertex is the one whose corresponding saddle connection has width equal to the sum of the other two widths. In the case of the $L^\infty$ Delaunay triangulation, we recover a train track that was dual to the Voronoi diagram, and such that the large edge at each vertex is dual to the saddle connection of largest width in each triangle.\\

\begin{definition}We refer to this train track as the \emph{vertical track of} $(X,q)$. For the same construction with the roles of vertical and horizontal interchanged we refer to this track as the \emph{horizontal track of} $(X,q)$. Similarly, we define vertical and horizontal tracks associated to any veering triangulation.\end{definition}

\noindent The tangential measure on this track is described in terms of period coordinates in the following manner:\\

\begin{definition}If $e$ is an edge of a train track, we say $f$ is a \emph{small partner} of $e$ if there is a vertex at which $e$ and $f$ are both small. Counting with multiplicity, an edge $e$ has zero, one, or two small partners if it is large, mixed, or small, respectively.\end{definition}

\begin{proposition}Let $(X,q)$ have a veering triangulation by saddle connections. Then the heights of the rectangles corresponding to edges in \Cref{PropCanonicalRectangles} are half the sum of the heights of their small partners, counted with multiplicity.\end{proposition}

\noindent Proof: We have already seen that large edges give rise to rectangles of height $0$. If an edge $e$ is small, then we may connect the midpoint of $e$ to the midpoint of each of its large neighbors by a segment of height equal to the small partner it meets at that large vertex; For either choice of up, the midpoint of one of the large neighbors will be above the midpoint of $e$ and the other will be below the midpoint of $e$. If $e$ is mixed, then the height of the segment joining the midpoint of $e$ and the midpoint of the large neighbor of $e$ is the height of the saddle connection dual to the small partner of $e$. $\Box$\\

\noindent Let $\tau$ be a train track with an associated transverse and tangential measure, $\mu$ and $\nu$, respectively, forming a quadratic differential. For each edge $i$ of $\tau$ let $R(i)$ be the associated rectangle. If $e$ is a large edge of $\tau$, let $a,b,d,c$ be its four neighbors starting at the top left and moving clockwise. If $$\mu(a) + \mu(d) > \mu(b) + \mu(c),$$ we may perform a \emph{left split} of the edge $e$ by moving the portions of the rectangle $R(e)$ that are vertically aligned with $R(b)$ and $R(c)$ into $R(b)$ and $R(c)$. The The result is a new train track and rectangle decomposition in which $e$ is small at both ends. $a$ and $d$ are said to be the \emph{winners}; $b$ and $c$ are the \emph{losers}. Similarly, if $$\mu(a)  +  \mu(d) < \mu(b) + \mu(c)$$ we may perform a \emph{right split}. The inverse operation of a split is a \emph{fold}. For a Teichm\"uller geodesic whose associated quadratic differentials are suited, the train tracks of \Cref{PropCanonicalRectangles} change by a left or right split at precisely those times when a maximal embedded square has four cone points on its boundary - one at each edge; the saddle connection that joined the left and right edges of the square is replaced by the saddle connection joining the top and bottom edges in the $L^\infty$ Delaunay triangulation.\\

\noindent We recall the following property of splitting \cite{PH}, Proposition 2.2.1:

\begin{proposition}If $\tau$ is birecurrent and $\tau^\prime$ is obtained from $\tau$ by performing a left or right split of some edge, and $\tau^\prime$ admits a strictly positive transverse measure, then $\tau^\prime$ is birecurrent.\end{proposition}

\begin{figure}[H]\label{splitpic}
\begin{tikzpicture}
\draw[ultra thick] (0,0) --(0,6);
\draw[ultra thick] (3,0) --(3,6);
\draw[ultra thick] (0,1) --(3,1);
\draw[ultra thick] (0,5) --(3,5);
\draw[ultra thick] (2,6) --(2,5);
\draw[ultra thick] (1,0) --(1,1);

\draw [dashed] (0.8,0) --(1,1) --(1.2,0);
\draw [dashed] (1.8,6) --(2,5) --(2.2,6);
\draw [dashed] (1,1) --(2,5);

\draw[ultra thick] (9,0) --(9,6);
\draw[ultra thick] (6,0) --(6,6);
\draw[ultra thick] (9,1) --(7,1);
\draw[ultra thick] (8,1) --(8,6);
\draw[ultra thick] (6,5) --(8,5);
\draw[ultra thick] (7,5) --(7,0);

\draw [dashed] (8.5,6) --(8,1) --(8.1,0);
\draw [dashed] (6.9,6) --(7,5) --(6.5,0);
\draw [dashed] (8,1) --(7,5);

\end{tikzpicture}

\caption{The second configuration of rectangles pictured is obtained from the first by performing a left split of the edge in the middle. The top left and bottom right rectangles are the winners of the split; the bottom left and top right are the losers. The associated train tracks are drawn dashed with nearly vertical edges.}
\end{figure}
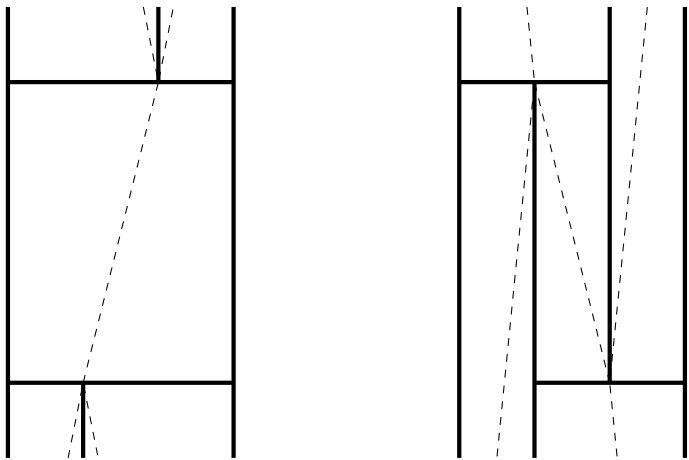

\noindent The following proposition concerns the evolution of veering triangulations by splitting of edges along Teichm\"uller flow lines. It is essentially due to Agol \cite{Veering} and Gu\'eritaud \cite{FG}. (Agol described the triangulation moves, and Gu\'eritaud described the surfaces on which the $L^\infty$ Delaunay triangulation is not unique.)

\begin{proposition}\label{Successor} If $(X,q)$ is suited, then along its Teichm\"uller flow line $g_t (X,q)$, the saddle connections that correspond to edges of train tracks that appear are precisely those that are diagonals of oriented rectangles. An saddle connection $e$ is replaced by a saddle connection $e^\prime$ at some point during the associated splitting sequence if and only if there is a rectangle whose vertical edges contain the endpoints of $e$ and whose horizontal edges contain the endpoints of $e^\prime.$ The corresponding split is a left split if and only if $e^\prime$ has positive slope, and it is a right split if and only if $e^\prime$ has negative slope. Moreover, if $\tau^\prime$ is obtained by $\tau$ by a split and they are the vertical tracks dual to veering triangulations $T$ and $T^\prime$, then [$e$ is an edge of $T$ but not of $T^\prime$] if and only if [$e^\prime$ is an edge of $T^\prime$ but not of $T$]. \end{proposition}

\noindent In \cite{VinCor}, Vincent Delecroix and Corinna Ulcigrai studied decompositions of hyperelliptic translation surfaces into quadrilaterals that evolve under the same types of edge replacement moves. Their purpose in doing this was to find an automaton (directed graph) that is much simpler than the veering automaton, i.e. the directed graph whose vertices are veering triangulations up to combinatorial equivalence, and an edge from $T_1$ to $T_2$ if $T_2$ is obtained from $T_1$ by the dual move to an edge split.\\

\noindent We remark that the condition that a quadratic differential have a particular train track is an open condition cut out by linear inequalities in the real and imaginary parts of the period coordinates. Throughout, while we have insisted on using suited quadratic differentials, the whole construction is valid if all vertical and horizontal saddle connections are sufficently long compared to the diameter of the surface. The space of quadratic differentials for which the construction of the rectangle decomposition of \Cref{PropVertSeparatrix} is therefore open, dense, and full measure with respect to a natural Lebesgue measure on the space of quadratic differentials, called the \emph{Masur-Veech measure}, which we describe now. In \cite{ietsandmfs} and \cite{Ve86}, Masur and Veech constructed this measure (in slightly different contexts) in order to apply ergodic theory to the study of the Teichm\"uller flow. If $(X,q)$ is a generic quadratic differential with cone point set $\Sigma$, the Masur-Veech measure is the natural Lebesgue measure on $\cc^n$ in period coordinates, provided that in the lifts of the saddle connections span the right subspace of the homology of an orienting double cover of $(X,q)$. This construction comes from \cite{ietsandmfs} and we will elaborate on it later; the important thing will be that almost every quadratic differential is suited. (There are only countably many saddle connections, so we can rotate by almost any angle to obtain a surface with no horizontal or vertical ones.)\\

\noindent Quadratic differential space $QD(\mathcal{T}_{g,n})$ is stratified based on the numbers and types of cone points: each \emph{stratum} is a maximal set where the number of cone points of each cone angle is constant, and such that for each $m \geq 1$, the number of marked points with cone angle $m\pi$ is constant. The \emph{principal stratum}, which consists of quadratic differentials with contains $n$
simple poles and $4g-4+n$ simple zeros, has full measure. Each stratum is an affine complex manifold, via the atlas of charts coming from systems of period coordinates.\\

\noindent We may summarize the construction of the $L^\infty$ Delaunay triangulation of $(X,q)$ and associated train track $\tau$ in the following way:

\begin{itemize}\item $\tau$ is the dual graph to the $L^\infty$ Delaunay triangulation.
\item Complementary $n$-gons of $\tau$ correspond to cone points of cone angle $n\pi$.
\item Vertices of $\tau$ correspond to vertical separatrices. Each vertex is a discontinuity in the inward normal vector of a unique region, which is the cone point whose separatrix it corresponds to.
\item If $q$ has $c$ cone points, $\tau$ has $2c + (4g-4)$ vertices and $3c + (6g-6)$ edges. This follows from the fact that $\tau$ is a 3-regular graph and the Euler characteristic of the surface is $2 - 2g$ if we do not delete the cone points.\end{itemize}

\section{A Greedy Algorithm for the $L^\infty$ Delaunay Triangulation}

\begin{definition} For a fixed stratum of quadratic differentials $\mathcal{Q}(\mu)$ let the \emph{veering triangulation automaton} associated to $\mu$ be the graph with the following: \begin{itemize}\item one vertex for each marked veering triangulation of a half-translation surface structure arising from quadratic differentials in $\mathcal{Q}(\mu)$ \item an edge for each pair of veering triangulations related by the following operation: Pick two triangles of the triangulation $T$ sharing an edge, such that the union of these two open triangles and the non-end points of their shared edge is an open quadrilateral whose boundary sides are edges of $T$ with alternating slope signs, and replace the diagonal of this quadrilateral that belongs to $T$ with the diagonal that does not belong to $T$. We will use $\mathrm{VT}(\mu)$ to denote this graph. For a fixed half-translation surface $(X,q)$ we let $\mathrm{VT}(X,q)$ denote the subgraph consisting of veering triangulations compatible with $(X,q)$\end{itemize}\end{definition}

\begin{theorem}\label{greedy}Suppose $(X,q)$ is suited and has a unique $L^\infty$ Delaunay triangulation. Then $\mathrm{VT}(X,q)$ is connected, and the $L^\infty$ Delaunay triangulation is characterized by the property that for each quadrilateral consisting of two triangles sharing an edge, such that the remaining four sides alternate in slope sign, the edge that the triangles share has shorter $L^\infty$ length than the quadrilateral's other diagonal. In particular, there is a greedy algorithm to reduce any veering triangulation to the $L^\infty$ Delaunay triangulation - find a quadrilateral for which this hypothesis is violated, and swap diagonals to reduce the sum of the lengths of the edges in the $L^\infty$ flat metric.\end{theorem}

\noindent $\mathrm{VT}(X,q)$ was previously known to be connected. Yair Minsky and Samuel Taylor proved this in section 3 of \cite{MinTay}, by a different argument.\\

\noindent Proof: First we show that the $L^\infty$ Delauany triangulation has this property. Let ABCD be a quadrilateral in the plane with sides alternating positive and negative slope, and let diagonal AC have shorter $L^\infty$ length than diagonal BD. Then the circumsquare of triangle BCD contains the point A. Therefore this pair of triangles could not have been a pair of neighboring triangles in the $L^\infty$ Delaunay triangulation of $(X,q)$.\\

\noindent Conversely, suppose that a veering triangulation has this property. Then, for each triangle in the universal cover, we must show that the circumsquare is embedded. Because the triangulation is veering, the vertices of each triangle are on the boundary of an embedded rectangle with vertical and horiztonal sides; we will call the minimal such rectangle the \emph{circumrectangle} of this triangle. We will show that the circumrectangle is contained in a square with vertical and horizontal sides and which has the the three vertices on the boundary; i.e. the circumrectangle is contained in the circumsquare, and therefore the triangle was an $L^\infty$ Delaunay triangle. Since this is true for all triangles we will conclude that the triangulation is the $L^\infty$ Delaunay triangulation.\\

\noindent Without loss of generality, assume that the width of the circumrectangle associated to some triangle ABC is greater than or equal to the height. (We lose no generality because we may interchange vertical and horizontal, and widths and heights, in the argument that follows.) Now, to the triangulation we have associated vertical and horizontal train tracks. Consider the vertical train track. There is a widest edge, say BC, of the triangle ABC, to which corresponds a large half-edge of a vertex. Follow a train route away from this vertex, starting along the large half-edge, for as long as possible without making any choices; this will end as soon as a large edge is reached (possibly the edge dual to BC itself, i.e. the first edge), and the weights on the edges only increase along the sequence of edges traversed. In particular, the vertical separatrix from A whose initial segment is contained in triangle ABC will pass through a sequence of increasing wider edges, and for each pair of consecutive edges, they will be two sides of a triangle of which the second edge is dual to a large half-edge of the vertical track. Eventually a large edge must be reached, since it is impossible to traverse the same edge twice in the same direction without passing through a large edge. By the diagonal property, the large edge is contained in a quadrilateral whose four edges lie on the boundary of a rectangle that is taller than it is wide. This implies that our original rectangle can be extended to a square.\\

\noindent To prove that the greedy algorithm works, we simply note that for any $M > 0$ on any $(X,q)$ there are only finitely many saddle that are of $L^\infty$ flat length less than $M$, so the algorithm must eventually terminate.$\Box$

\begin{corollary}Given a veering triangulation and an assignment of heights and widths, a finite collection of inequalities in the heights and widths certify whether or not the triangulation is the $L^\infty$ Delaunay triangulation.\end{corollary}

\noindent Remark: It is possible to write a finite orbifold cover of $QD(\mathcal{M}_{g,n})$ as a finite union of convex sets with disjoint interiors, based on the topological type of the $L^\infty$ Delaunay triangulation. However, already in the case $\mathcal{M}_{1,1}$, which has lowest complexity of all moduli spaces one might try to consider, a cell of top dimension may have infinitely many cells of lower dimension on its boundary. (The author is grateful to Simion Filip for pointing out this complication.) However, one might expect that the infinite families can be parametrized by finite data, for instance; one may hope that all neighbors of a fixed cell are obtained from one of finitely many other cells up to Dehn twists.

\section{Short Curves on Riemann Surfaces}

\noindent A simple closed curve $c$ in $S_{g} \setminus \{x_1,...,x_n\}$ is \emph{essential} if neither component of $S_g \setminus \{x_1,...,x_n\} \setminus c$ is homeomorphic to a disk or once-punctured disk. We will refer to isotopy classes of unoriented essential simple closed curves as \emph{curves}. The \emph{curve complex} of $S_{g,n}$ is a simplicial complex $\mathcal{C}(S_{g,n})$ whose vertices are curves, and whose simplices correspond to sets of curves that can be realized by nonintersecting simple closed curves. We can metrize the curve complex by making each simplex the standard simplex in $\rr^d$ with edge length $1$. In \cite{MM}, Masur and Minsky showed that the curve complex of $S_{g,n}$ is Gromov hyperbolic whenever $3g-3+n > 1$.\\

\noindent Let $\tau$ be a filling birecurrent train track admitting transverse and tangential measures and an associated rectangle decomposition $R(\tau)$. A \emph{multicurve carried by} $\tau$ is an isotopy class of finite union of simple closed curves that do not intersect the vertical edges of the rectangles $R(\tau)$, and such that each component of an intersection of each curve with a rectangle $R(i)$ meets $\partial R(i)$ at exactly two points, one on each horizontal edge, and intersects those edges transversely. A multicurve is equivalent to an integer-valued transverse measure on $\tau$, except that the weights of some edges can be allowed to be zero (but not all of them.) We do not require that the components of a multicurve belong to different isotopy classes of simple closed curve. A multicurve is equivalent to a transverse measure that assigns non-negative integer weight to each edge (but possibly zero weight to some edges). The transverse measures with total edge sum $1$ form a convex polytope whose extremal points are rational points, because the polytope is cut out by rational linear equations. The minimal integral multiples of the vertices of this polytope are called \emph{vertex curves} of $\tau$. A vertex curve is determined by a periodic sequence of oriented edges of $\tau$, and in each period it visits each edge at most twice, with different orientations if it visits twice.\\

\noindent We have the following special case of Rafi's thick-thin decomposition, allowing us to compare the hyperbolic and flat lengths of geodesics. If $\alpha$ is a curve, $(X,q) \in QD(\mathcal{T}_{g,n})$, let $\ell_q(\alpha)$ be the infimal length of curves in the isotopy class $\alpha$ in the $q$-metric, and let $\ell_{\sigma_X}$ be the length of the geodesic in the isotopy class of $\alpha$ in the hyperbolic metric. We have the following theorem from \cite{ThickThin}.

\begin{theorem}\label{CompareLength} Fix nonnegative integers $g,n$ with $3g-3+n < 0$. Fix $\epsilon > 0$. Then there is a constant $C_\epsilon \in [1,\infty)$, with the following property: if $(X,q) \in QD^1(\mathcal{T}_{g,n})$ be such that every closed geodesic $\gamma$ has $\ell_q(\gamma) \geq \epsilon$, then for all curves $\alpha \in \mathcal{C}(S_{g,n})$ we have $$C_\epsilon^{-1} \ell_q(\alpha) \leq \ell_{\sigma_X}(\alpha) \leq C_\epsilon \ell_q(\alpha).$$ \end{theorem}

\noindent By Mumford's compactness criterion, \cite{Mumf}, the subset of $\mathcal{M}_{g,n}^{(\epsilon)} \subsetneq \mathcal{M}_{g,n}$, the set of Riemann surfaces whose associated hyperbolic metrics have no curves of length $< \epsilon$, is compact. In particular, this means that we may pick compact set $K \subset \mathcal{T}_{g,n}$ whose projection contains every element of $\mathcal{M}_{g,n}^{(\epsilon)}$. If $K$ has diameter $M$, then for any $X,Y \in K$ we have $$e^{-2M} \leq \frac{\ell_{\sigma_X}(\alpha)}{\ell_{\sigma_Y}(\alpha)} \leq e^{2M}$$ for all curves $\alpha$, by Wolpert's lemma \cite{Wolp}. For each fixed $X$ there are only finitely many curves of length $D$, even if $X$ has cusps, since all essential simple closed curves meet a compact subset of $X$ and the deck group of $X$ acts properly discontinuously.\\

\noindent We thus have the following:

\begin{lemma}\label{DelaunayMarking} Fix $\epsilon > 0$. Then there is a constant $C_\epsilon < \infty$ such that whenever $(X,q) \in QD(\mathcal{M}_{g,n}^{(\epsilon)})$ the curve complex distance between the shortest curve in the metric $\sigma_X$ and each vertex curve of the train track $\tau_q$ associated to the $q$ is at most $C_\epsilon$. \end{lemma}

\noindent Proof: since the Teichm\"uller metric is proper, we may pick a closed bounded ball in $\mathcal{T}_{g,n}$ whose projection covers $\mathcal{M}_{g,n}^{(\epsilon)}$. As $(X,q)$ varies in this ball, the diameter of $(X,q)$ remains bounded. The diameter of an $L^\infty$ Delaunay triangle in $(X,q)$ is bounded by $2\sqrt{2}$ times the diameter of the surface. The train track and its complementary regions is Poincar\'e dual to the $L^\infty$ Delaunay triangulation, and each vertex curve visits each edge of $\tau_q$ at most twice, and the number of edges depends only on $g$ and $n$, this implies a uniform bound on the $q$-lengths of vertex curves. By \Cref{CompareLength} and local finiteness of curves of hyperbolic length $D$, the lemma follows.\\

\noindent Let $S_{g,n}$ be a surface of genus $g$ with $n$ marked points $\{x_1,...x_n\}$. The \emph{mapping class group} $\mathrm{Mod}(S_{g,n})$ is the quotient group of orientation preserving homeomorphisms of $S_{g,n} \setminus \{x_1,...,x_n\}$ that fix the $n$ marked points by the subgroup homotopic to the identity. Elements of $\mathrm{Mod}(S_{g,n})$ are called \emph{mapping classes}. The mapping class group acts properly discontinuously on $\mathcal{T}_{g,n}$ with quotient $\mathcal{M}_{g,n}$, see e.g. \cite{FarbMarg}.

\begin{definition} A \emph{pseudo-Anosov map} from $S_{g,n}$ to itself is a homeomorphism with the following properties: \begin{itemize} \item \noindent $f$ maps marked points to marked points \item $S_{g,n}$ admits a half-translation surface structure, whose only cone points of angle $\pi$ occur at marked points, with respect to charts on which $\sqrt{q} = \pm(dx + i dy)$, the derivative of $f$ at every non-cone point is of the form $$g_t := \pm\left(\begin{array}{rl} e^t & 0\\ 0 & e^{-t} \end{array}\right) \in \psltwo$$ for some $t > 0.$ \end{itemize} A mapping class containing a pseudo-Anosov map is called a \emph{pseudo-Anosov mapping class}.\end{definition}

\noindent Pseudo-Anosov maps give rise to closed geodesics in $\mathcal{M}_{g,n}$.\\

\noindent In particular, a pseudo-Anosov mapping class is never realized by a homeomorphism supported on a subsurface whose boundary contains an essential curve.

\section{Measured Foliations, Positive Cones, and Hilbert Metrics}

\noindent Previously we defined spaces of transverse and tangential measures for a complete birecurrent train track. A transverse measure together with an equivalence class of tangential measure induces a quadratic differential on a Riemann surface, which comes equipped with vertical and horizontal \emph{measured foliations}. Measured foliations arise from quadratic differentials as singular foliations of the surface $S_g$ by vertical and horizontal leaves, with a singularity at each cone point; more precisely, a measured foliation is equivalent to a maximal system of charts with transition maps that preserve a $1$-form up to sign $|dx|$, and the $1$-form is locally smooth and non-vanishing at all but finitely many points, called \emph{singularities} of the measured foliation; it may be undefined at these points. We also require further that at each point where the 1-form is not defined up to sign, the foliation must have a neighborhood isotopic to the quotient of the foliation of $\cc$ given by level sets of the function $Re(z^k)$ by the map identifying $z$ with $-z$, which preserves the $1$-form up to sign away from $0$; a singularity of this foliation is said to be $k$-pronged. A $2$-pronged singularity is in some sense not a singularity; if these occur, it is because the foliation extends smoothly but the point is a marked point. A 1-pronged singularity is permitted at each puncture/marked point, and at all other points the foliation is smooth or has a $k$-pronged singularity for some $k > 2$. For a quadratic differential $(X,q)$, we say that $Re \sqrt{q}$ is the associated 1-form up to sign for the \emph{vertical foliation} of $(X,q)$, and $Im \sqrt{q}$ is the 1-form for the \emph{horizontal foliation}.\\

\noindent Two measured foliations are \emph{Whitehead equivalent} if there is a third measured foliation that can be obtained from both of them, up to isotopy, by collapsing compact segments of leaves whose endpoints are singularities, at most one of which is a marked point. If $\tau$ is a complete birecurrent train track, and $(\mu_1,\nu_1)$ and $(\mu_2, \nu_2)$ are ordered pairs consisting of a transverse measure in the first coordinate and a tangential measure in the second, then the associated quadratic differentials $(X_1,q_1)$ and $(X_2,q_2)$. $\mu = \mu^\prime$ if and only if the vertical foliations are Whitehead equivalent; similarly, $\nu_1$ and $\nu_2$ belong to the same equivalence class of tangential measure if and only if the horizontal foliations of the two quadratic differentials are Whitehead equivalent.\\

A \emph{strongly stable leaf} for the Teichm\"uller geodesic flow is a maximal set of unit area quadratic differentials whose vertical measured foliations belong to one Whitehead equivalence class. If $\tau$ is a complete birecurrent train track, the set $\mathcal{MF}$ of all measured foliations carried by $\tau$ is a manifold of dimension $6g-6+2n$. We will describe a system of charts for $\mathcal{MF}$ after we state the Hubbard-Masur theorem, \Cref{HubMas}.\\

A {stable leaf} for the Teichm\"uller flow is a maximal set of unit area quadratic differentials whose vertical measured foliations are equivalent up to scaling the measure and Whitehead moves. The \emph{strongly unstable} and \emph{unstable leaves} are defined in the same way, with the horizontal foliation in place of the vertical foliation.\\

\noindent The following theorem of Hubbard and Masur \cite{HM} tells us how the leaves of these foliations sit inside of Teichm\"uller space:

\begin{theorem}\label{HubMas} Let $X \in \mathcal{T}_{g,n}$ be a Riemann surface. Each Whitehead equivalence class of measured foliation arises from a unique nonzero quadratic differential on $X$. Moreover, $QD(\mathcal{T}_{g,n})$ is, as the cotangent bundle of $\mathcal{T}_{g,n}$ (minus the zero section), naturally homeomorphic to $\mathcal{T}_{g,n} \times \mathcal{MF}$ via the map that assigns to each quadratic differential its underlying Riemann surface and vertical foliation.\end{theorem}

\noindent The topology on $\mathcal{MF}$ in the Hubbard-Masur theorem (\Cref{HubMas}) is the following: the space of transverse measures (weights satisfying the switch condition) on complete birecurrent train tracks is a chart on the space $\mathcal{MF}$ to $\rr^{6g-6+2n}$, and the collection of all such charts covers $\mathcal{MF}$ in such a way that the transition maps between charts are all linear maps. There is a natural multiplication action by $\rr^+$ which is just to multiply the measure by positive scalars; the quotient by this action is the space of \emph{projective measured foliations}, or $\mathcal{PMF}$, and it is homeomorphic to the sphere $S^{6g-7+2n}$. The intersection number of a curve $c$ and a measured foliation given by a 1-form up to sign $|dx|$ is given by $$i(c,dx) = \inf\limits_{[\gamma] = c} \int_\gamma |dx|,$$ and every measured foliation is determined by its intersection numbers with a finite collection of curves. If the underlying surface has marked points or punctures, two simple closed curves do not represent the same curve unless they are connected by an isotopy through curves that do not pass through any marked points/punctures. The mapping class group acts on $\mathcal{MF}$ by linear homeomorphisms with respect to train track coordinates, and acts minimally on the quotient $\mathcal{PMF}$ by homeomorphisms. See \cite{FLP} for precise definitions and classification.

\begin{definition}A \emph{positive cone} is a closed subset $C$ of a finite dimensional real vector space satisfying all of the following conditions: \begin{itemize}
    \item $C \cap -C = \{0\}$ \item $\lambda C = C$ for all $\lambda > 0$ \item $C$ is closed under vector addition.
\end{itemize}\end{definition}

\begin{proposition} Let $\tau$ be the vertical train track associated to a veering triangulation of a quadratic differential in the principal stratum, and let $\rr^{E(\tau)}$ the vector space of functions on the edges of $\tau$, with standard basis $$\{\mathbb{1}_e: e \mathrm{~is~an~edge~of~} \tau\}.$$ Let $\rr^{V(\tau)}$ be the space of differences between equivalent transverse measures (recall that this space is generated by vectors of the form $\mathbb{1}_a - \mathbb{1}_b - \mathbb{1}_c$ where $a$ is the large edge at a vertex with small edges $b$ and $c$). The following are positive cones: \begin{itemize}\item The closure of the space of tangential measures on $\tau$ in $\rr^{E(\tau)}$.\item The closure of the space of transverse measures carried by $\tau$ in $\rr^{E(\tau)}$ \item The projection of the closure of the space of tangential measures to the quotient of $\rr^{|E(\tau)|}$ by $V(\tau)$.\end{itemize} \end{proposition}

\noindent Proof: The space of transverse measures forms a positive cone because it is the intersection of the space of non-negative functions on the edges with a linear subspace. The closure of the space of tangential measures is a positive cone because it is an intersection of the space of non-negative functions on the edges with a collection of half spaces (given by three linear inequalities for each complementary $3$-gon). So all we need to check is the quotient of the tangential measures by the space of differences of equivalent tangential measures. The only condition that is not obvious is the first condition, $C \cap -C = \{0\}$.\\

\noindent For this, we consider the standard Euclidean inner product on $\rr^{E(\tau)}$, and assume that the coordinates correspond to edges. Transverse and tangential measures are functions on the edges. A transverse measure and a tangential measure have inner product equal to the area of the associated quadratic differential. If $\mu$ is any transverse measure with positive weight on every edge, which must exist, and $\nu$ is any nonzero nonnegative function, then $(\mu,\nu) > 0$, and $(\mu,\nu^\prime) = 0$ if $\nu^\prime$ is a difference of equivalent tangential measures. This linear functional on $\rr^{E(\tau)}/\rr^{V(\tau)}$ implies the desired property. $\Box$\\

\begin{notation}Let $\tau$ be a complete birecurrent train track. Write $\mathcal{MF}(\tau)$ for the cone of equivalence classes of tangential measures on $\tau$, and $\mathcal{MF}^\pitchfork(\tau)$ for the cone of transverse measures on $\tau$.\end{notation}

\noindent We remark that by the hyperplane separation theorems, any positive cone $C$ is separated from $-C$ by a hypeplane, i.e. there is a linear functional that is positive on $C \setminus \{0\}$, and the set $D_\Lambda := \{x \in C: \Lambda x = 1\}$ is convex and homeomorphic to a closed disk.\\

We recall the definition of Hilbert metric and its properties with respect to maps between cones. These are summarized by Birkhoff in \cite{Birk}.\\

\begin{definition}Let $C$ be a positive cone, and let $\mathbb{P}(C)$ be the space of rays from $0$ contained in $C$. The \emph{Hilbert pseudo-metric} on $\mathbb{P}(C)$ is defined in the following manner: If $r_1$ and $r_2$ are two distinct rays emanating from $0$ in $\mathbb{P}(C)$, and $\ell$ is a line parametrized proportional to arc length that intersects $C$ precisely on the interval $[0,1]$, meeting $r_1$ at $a \in [0,1]$ and $r_2$ at $b \in [0,1],$ then $$d_H(r_1,r_2) := \left|\log\frac{a(1-b)}{b(1-a)}\right|.$$ \end{definition}

Since this is the log of a cross ratio of four points in the projectivization of the plane spanned by $r_1$ and $r_2$, it does not depend on $\ell$. $\mathbb{P}(C)$ is naturally a manifold with boundary (in fact a manifold with corners); it can be identified with the level set of some functional that is positive on $C$. The Hilbert pseudo-metric is a metric on the interior of $\mathbb{P}(C).$ Two points on the boundary of $\mathbb{P}(C)$ may be finite or infinite distance apart, but the distance from an interior point to a boundary point is infinite.\\

\noindent The key properties of Hilbert metrics are the following:\\

\begin{theorem}\cite{Birk} Let $f:\rr^{m_1} \to \rr^{m_2}$ be a linear map, and suppose $C_i$ is a positive cone in $\rr^{m_i}$ with associated Hilbert metric $d_{H_i}$. If $f(C_1) \subset C_2$, then for all rays $x,y \in \mathbb{P}(C_1),$
$$d_{H_2}(f(x),f(y)) \leq d_{H_1}(x,y).$$ Moreover, if $f(\mathbb{P}(C_1))$ has finite diameter with respect to $d_{H_2}$, then $f$ is a contraction mapping, i.e. there is a number $c < 1$ such that for all $x,y \in \mathbb{P}(C_1)$, $$d_{H_2}(f(x),f(y)) < c d_{H_1}(x,y).$$

\noindent If $\Delta$ is the diameter of $f(\mathbb{P}(C_1))$ we may take $c = \tanh(\Delta/4).$
\end{theorem}

\noindent We add in one extra observation, which will be useful later.\\

\begin{corollary}\label{CoroGoodWord}
Suppose that a linear map $T:\rr^m \to \rr^m$ maps a positive cone $C$ into itself, and $T^k$ maps every ray of $\mathbb{P}(C)$ into the interior of $\mathbb{P}(C)$. Let $S_1,...,S_{k-1}: \rr^m \to \rr^m$ be linear maps with the property that for all $x \in C$, $S_i(x) - x \in C$. Then $T \circ \prod\limits_{i = 1}^{k - 1} S_i \circ T$ is a contraction with respect to the Hilbert metric $d_H$.\\
\end{corollary}

\noindent Proof: An easy induction on $k$ shows that for all $x \in C$, $\left(T \circ \prod\limits_{i = 1}^{k - 1} S_i \circ T\right)x - T^k x \in C.$ Since $C$ is the cone on a closed disk $D_\Lambda$ and every point of $D_\Lambda$ maps to an interior point of $C$, it follows that the disk $D_\Lambda$ maps onto a compact subset of the interior of $C$. The metric $d_H$ is bounded on this set so we're done. $\Box$\\

\section{The Transition Matrices of a Splitting Sequence}

In this section, we consider a periodic train track splitting sequence. Let $\{\tau_m\}_{m = -\infty}^\infty$ be a bi-infinite sequence of \emph{labelled} train tracks, i.e. the edges, vertices, and complementary regions have names, and such that $\tau_{m+1}$ is obtained from $\tau_m$ by splitting one or more large edges. In terms of rectangle decompositions, this means that rectangles, cone points, and separatrices have names.\\

\noindent In terms of an associated rectangle decomposition coming from transverse measure $\mu_m$ and tangential $\nu_m$, $\tau_{m+1}$ inherits a rectangle decomposition from the splitting, by increasing heights of rectangles associated to the labels of losing edges, decreasing the width of the rectangle associated to the label of the splitting edge, and keeping all other rectangles the same. (See \Cref{splitpic}.) Then we have the following:

\begin{proposition}\label{BasicSplitMatrix} Let $TGM(\tau)$ denote the cone of tangential measures on $\tau$. \emph{Note that we are not taking equivalence classes of tangential measures, hence the different notation.} If $\tau_{m+1}$ can be obtained from $\tau_{m}$ by splitting an edge of $\tau_m$ for some sequence of train tracks $\{\tau_m\}_{m = r}^s$, $-\infty < r < s < \infty.$ Let $E$ be an index set for the edges of any track in the sequence (and therefore for all tracks in the sequence). Then there are induced linear maps on $\rr^E$ which induce maps $\phi_{m,m+1}^{\pitchfork}: \mathcal{MF}^\pitchfork(\tau_{m+1}) \to \mathcal{MF}^\pitchfork(\tau_m)$ and $\phi_{m,m+1}^{TG}: TGM(\tau_m) \to TGM(\tau_{m+1})$. These linear maps on $\rr^E$ have the following description:\\

\noindent Let $I$ be the identity matrix on $\rr^E$, and let $M_{ij}$ be the elementary matrix whose unique nonzero entry is a $1$ in entry $ij$. With respect to the standard basis for $\rr^E$, the maps $\phi_{m,m+1}^{\pitchfork}$ and $\phi_{m,m+1}^{TG}$ are induced by the following matrices: For each edge $e$ that splits with winners $a,d$ and losers $b,c$, we have
$$\phi_{m,m+1}^{\pitchfork} = I + M_{eb} + M_{ec},\phi_{m,m+1}^{TG} = I + M_{be} + M_{ce}.$$\end{proposition}

\noindent Proof: Obvious from the definition of splitting. $\Box$\\

\noindent Two important corollaries follow:

\begin{corollary}\label{Monotone}If $m_1 < m_2$, then the sequence $\tau_m$ admits a pair of invertible matrices $\phi_{m_1,m_2}^\pitchfork: \mathcal{MF}^\pitchfork(\tau_{m_2}) \to \mathcal{MF}^\pitchfork(\tau_{m_1})$ and $\phi_{m_1,m_2}^{TG}: TGM(\tau_{m_1}) \to TGM(\tau_{m_2}).$ The matrices $\phi_{m_1,m_2}^\pitchfork - I$ and $\phi_{m_1,m_2}^{TG} - I$ have non-negative integer entries, and $\phi_{m_1,m_2}^\pitchfork$ and $\phi_{m_1,m_2}^{TG}$ induce the maps of transverse and tangential measures associated to the splitting sequence $\{\tau_m\}$.\end{corollary}

\noindent Proof: Induction on $m_1 - m_2$.$\Box$

\begin{corollary}\label{UniquelySplit} Let $\tau_{m_2}$ be obtained from $\tau_{m_1}$ by a sequence of splits, and assume both are generic and birecurrent. Then the composition of induced maps on transverse measures and induced maps on tangential measures are induced by a pair of unimodular integer matrices, which are adjoints. Moreover, these matrices can be recovered from a smaller collection of data: the initial track $\tau_1$ with an indexing of the edges, the final track $\tau_2$ with the indexing inherited from the sequence, and for each edge index, the word in the letters $L$ and $R$ corresponding to the sequence of left and right splits that that edges corresponding to that index perform along the splitting sequence.\end{corollary}

\noindent Proof: The only claim that is not obvious is the last one, since the matrices coming from single edge-splits are transposes of each other and are composed in the opposite orders. Now, to prove the last claim, we consider the following directed graph $\Gamma$ on the set of splits that occur: there is an edge from a splitting of $e$ to the next splitting $f$ if $f$ is a winner when $e$ splits. The reason for doing this is that $f$ cannot become large until $e$ splits, so the split of $e$ must occur before the split of $f$. Thus there is a partial order on the splits, based on the existence of a directed path from one split to another in $\Gamma$. We say that a split $s$ is \emph{ready to happen} if it is minimal in this partial order, i.e. a directed path in $\Gamma$ cannot start at another split and reach $s$.\\

\noindent We claim that possible splitting sequences from $\tau_1$ to $\tau_2$ that follow our splitting data simply correspond to repeatedly performing the following operation: choose a split that is ready to happen, perform it, and then remove the corresponding vertex from $\Gamma$. We will actually show how to recover the matrix $TV^{-1}.$ We begin by assigning a formal variable $x_i$ that represents the initial transverse measure of edge $e_i$, and we do not impose any relation between the values $e_i$. In particular, we ignore the switch condition for the sake of computing the matrix. For weights satisfying the switch condition, we have some choice in how we write the weight of an edge after a split as a function of the weight before, but in order to keep consistent with our transition matrices, we insist that we apply the inverse matrix of the basic transverse measure splitting matrix from \Cref{Monotone} at each step. This is the matrix that subtracts the weights of the losers from the weight of the splitting edge.\\

\noindent We claim that when $s$ is ready to happen, the following all depend only on the splits from which there are paths to $s$: the labels of neighbors of the edge $e$ splitting at $s$, including which pairs meet at each vertex, the clockwise order of the four vertices, the weights on $e$ and its four neighbors when $s$ is ready to happen as linear functions of the $x_i$. The proof is just induction on the number of splits that occur before $s$ is ready to happen. Indeed, the only time the weight of an edge changes is when it splits, and the only time the neighbors of an edge change is when it is one of the five edges involved in a split. Now, consider the most recent split affecting each end of $e$ before $s$. If an end of $e$ never split, $e$ must have always been large at that end, and the edges at that end of $e$ simply have their initial weights and labels. If there was a previous split $s^\prime$, then $e$ was the winner, then the inductive hypothesis applies to $s^\prime$, and the weights of the three edges at that vertex before and after the splitting $s^\prime$ are determined as an integral combination of the $x_i$. The values of the final weights as a function of the $x_i$ are therefore uniquely determined, since the weight of an edge only changes when it splits, and the weight of each edge is uniquely determined, as a linear function of the $x_i$. Uniqueness of $TV^{-1}$ follows.\\

\begin{proposition}\label{PropPresTriIneq}The triangle inequality for tangential measures is preserved by transition matrices for splitting sequences. In fact, if $h_1,h_2,h_3$ are heights of three sides of a complementary triangle for a train track with tangential measure, then $h_1 + h_2 - h_3$ can only increase along a split.\end{proposition}

\noindent Proof: a split adds the same quantity to two of the three side lengths and leaves the third unchanged. $\Box$

\begin{corollary}\label{CoroGoodWordPrep} Suppose a splitting sequence takes a labelled train track to another train track in the same orbit of the mapping class group. If the cones of tangential measures are identified to a cone $C$ via the homeomorphism, then the image each point $p$ in the induced map on tangential measures is the sum of $p$ and another point in $C$. \end{corollary}

\noindent Proof: This follows from \Cref{PropPresTriIneq} and \Cref{Monotone}.

\section{Resolving Collisions of Singularities}

\noindent In this section we discuss collisions of singularities and their relation to the $L^\infty$ Delaunay triangulations of nearby quadratic differentials in the principal stratum. Recall that complementary regions of a train track $\tau$ correspond to zeros and poles of a quadratic differential $q$; in what follows, our focus will be on collections of edges that can be deleted to collapse a small cluster of singularities into a single singularity.\\

\noindent Suppose $\{(X_m,q_m)\}$ is a sequence of quadratic differentials in the principal stratum of $QD(\mathcal{T}_{g,n)}$ for which the $L^\infty$ Delaunay triangulation is defined, and they converge to a quadratic differential $(X_\infty,q_\infty)$ that does not belong to the principal stratum, but which is suited. In particular, passing to a subsequence we may assume that $(X_m,q_m)$ all have isomorphic $L^\infty$ Delaunay triangulations, isomorphic not only as CW-structures but also such that all of the isomorphisms send positively sloped saddle connections to positively sloped saddle connections, and such that the period of each saddle connection converges. Also assume that the maps $f_m$ to the base surface push these triangulations forward to isotopic triangulations of the base surface $S_{g,n}$, which are isotopic via isotopies that fix the marked point as vertices. Let $\tau$ be a marked labelled train track associated to these isomorphic $L^\infty$ Delaunay triangulations. Pick $\epsilon > 0$ such that the ball of radius $2\epsilon$ about each singularity of $(X_\infty,q_\infty)$ is a disk containing no other cone points. We will refer to all saddle connections belonging to the $L^\infty$ Delaunay triangulation of the $\{q_m\}$ and having length at most $\epsilon$ in all but finitely many $(X_m,q_m)$ as \emph{vanishing saddle connections}.

\begin{proposition}\label{DyingVertices} For each $L^\infty$ Delaunay triangle of $(X_m,q_m)$, the large half-edge of the corresponding vertex of $\tau$ is dual to a vanishing saddle connection if and only if the two small edges are also dual to vanishing saddle connections.\end{proposition}

\noindent Proof: Let $a,b,c$ be the sides of triangle, and assume that $c$ is large, in which case we have $$\left|Re\left(\int_a \sqrt{q}\right)\right| + \left|Re\left(\int_b \sqrt{q}\right)\right| = \left|Re\left(\int_c \sqrt{q}\right)\right|.$$

Then if $\left|Re\left(\int_a \sqrt{q}\right)\right| \to 0$, $\left|Re\left(\int_b \sqrt{q}\right)\right| \to  0$ and $\left|Re\left(\int_c \sqrt{q}\right)\right| \to 0$. This would imply $q_\infty$ has a vertical saddle connection if $a,b$ are not both vanishing. This violates our assumptions. The converse follows for similar reasons. $\Box$

\begin{proposition}\label{CarryNothing} Let $\alpha$ be a simple closed curve that does not pass through any cone points. Assume $\alpha$ does not intersect the $L^\infty$ Delaunay triangulation except for finitely many transverse intersections with vanishing saddle connections. Then $\alpha$ is not essential. \end{proposition}

\noindent Note that the homotopy class of $\alpha$ only depends on the sequence of intersections.\\

\noindent Proof: For any $\epsilon > 0$, for all sufficiently large $m$ we can find a representative for the isotopy class $[\alpha]$ on $(X_m,q_m)$ that stays within distance $\epsilon$ of the union of the cone points. However, as $(X_m,q_m) \to (X,q)$ the infimum of the lengths of essential simple closed curves in the $q_m$ metrics converges to the infimal length of simple closed curves in the $(X_\infty,q_\infty)$ metric, which is nonzero. In particular, the diameter of the $\alpha$ in the $q_m$ metric goes to $0$, whence the length of $\alpha$ in the $q_m$ metric goes to $0$. By Rafi's estimate \Cref{CompareLength}, the hyperbolic length of $\alpha$ goes to $0$ in the hyperbolic metrics $\Sigma_{X_m}$. Thus by Wolpert's lemma, $\alpha$ cannot have positive hyperbolic length on $(X_\infty,q_\infty)$, and so $\alpha$ is not essential.\\

\noindent These properties of vanishing saddle connections motivate the following definition:

\begin{definition}Let $\tau$ be a complete birecurrent train track. A subgraph $H$ of $\tau$ is \emph{inessential} if every simple closed curve that is homotopic to a curve contained in $H$ is non-essential, and if for each vertex, the large half-edge belongs to $H$ only if both small half-edges do. Equivalently, $H$ is inessential if after deleting all edges of $H$, deleting all vertices of degree $0$ that remain, and merging the two edges at each valence $2$ vertex, the result is a filling birecurrent train track. We call this birecurrent train track the \emph{complement} of $H$ and denote it by $\tau - H$.\end{definition}

\begin{proposition}If $\tau$ is a complete birecurrent train track, and $H$ is inessential, and $\tau^\prime$ is obtained from $\tau$ by splitting an edge of $H$, and the corresponding bijection of edges sends $H$ to $H^\prime$ itself, then $H^\prime$ is also inessential. The same is true if an edge that does not belong to $H$ splits and at least one loser is in $H$.\end{proposition}

\noindent Proof: $H$ and $H^\prime$ have the same complement.$\Box$

\begin{proposition} Given an inessential subgraph $H$ of a birecurrent track $\tau$, $\tau$ can be reduced to the complement of $H$ by deleting one small edge at a time.\end{proposition}

\noindent Proof: A train route in $H$ cannot traverse an edge twice in the same direction, since if it did, then some curve carried by $\tau$ would be supported in $H$. However, by  and a birecurrent track never carries a simple closed curve that is not essential. (See \cite{PH}). This contradicts \Cref{CarryNothing}. Consider a maximal train route $\gamma$ contained in an inessential sugraph $H$. Since both ends of $H$ are small half-edges, and whenever $\gamma$ changes between two edges, it moves from a large half-edge to a small half-edge or vice-versa, there are 2 more small half-edges than large half-edges in $\gamma$. Thus $H$ contains a small edge, which can be removed. Now we may induct.$\Box$

\begin{proposition}\label{ShrinkToAxis} If $\{(X_m,q_m)\} \to (X_\infty,q_\infty)$ is a convergent sequence satisfying the assumptions of \Cref{DyingVertices}, then the $L^\infty$ Delaunay triangulation of $q_\infty$ is obtained from the $L^\infty$ Delaunay triangulation of the $q_m$ by the following combinatorial construction: Delete the interior of any $L^\infty$ Delaunay triangle with at least one edge that is a vanishing saddle connection. Collapse each vanishing saddle connection to a point, and identify pairs of non-vanishing saddle connections that form two out of three sides of an $L^\infty$ Delaunay triangle whose interior was deleted by an affine homeomorphism.\end{proposition}

\noindent Proof: We will assume each $X_m$ is compact with $n$ marked points, and let $\{(\tilde{X}_m,\tilde{q}_m)\}$ be metric normal pole-free-disk covers of the homeomorphic surfaces $X_m$, with one of the marked points fixed as a base point. On these covers, the distance between pairs of lifts of marked points in $X_m$ stays bounded below. The identifications of the $L^\infty$ Delaunay triangulations and the choice of base point fix a simplicial homeomorphism between any $X_{m_1}$ and $X_{m_2}$ for any $m_1,m_2$; we will think of $\{(\tilde{X}_m,\tilde{q}_m)\}$, equipped with the $L^\infty$ flat metric, as a sequence of metrics on a fixed simplicial complex, which has a fixed deck group action, and let $(\hat{X},\hat{q})$ be the lift of the quotient space described. $(\hat{X},\hat{q})$ naturally has a triangulation and a half-translation surface structure, since in the limit, each pair of identified edges has the same length and slope, and it has an action by $\pi_1(X)$ whose quotient is $(X_\infty,q_\infty)$. We will check that the $L^\infty$ Delaunay triangles of $\tilde{X}$ map to $L^\infty$ triangles of $\hat{X}$ or lower dimensional cells. If we do that, we are done, since the number of triangles in the cell complex $\hat{X}/\pi_1(X)$ is equal to the number of triangles in the $L^\infty$ Delaunay triangulation of $(X_\infty,q_\infty).$\\

\noindent The key observation is the following: $L^\infty$ flat metrics associated to the $q_m$ are given by the max of two sequences of pseudo-metrics that converge uniformly on compact sets. Therefore the metrics form a convergent sequence of functions on $\tilde{X} \times \tilde{X}$ (the metrics extend to the completion), uniformly convergent on compact sets, and their limit is the pseudo-metric obtained by pulling back the $\hat{q}_\infty$ $L^\infty$ flat metric on $\hat{X}$ by the quotient map. The locations of the circumcircles form a Cauchy sequence with respect to this pseudo-metric on $\tilde{X}$, as do distances to all cone points. In particular, this implies that there are three cone points equidistant to the circumcenter of each circumsquare and no cone point is closer than these three. Since $q_\infty$ is suited, this implies the triangles are $L^\infty$ Delaunay triangles. $\Box$\\

\noindent By a similar argument, we have the following:

\begin{proposition}\label{DegenCell} Let $(X_m,q_m)$ be a convergent sequence of quadratic differentials sharing a veering triangulation $T$. Suppose $(X_m,q_m) \to (X_\infty,q_\infty)$ and $(X_\infty,q_\infty)$ has no horizontal or vertical saddle connections. Then the periods of all saddle connections in $T$ converge to periods of saddle connections in $(X_\infty,q_\infty)$ or to $0$ and the saddle connections whose periods go to $0$ are the edges of an inessential subgraph. Moreover, $X_\infty,q_\infty$ inherits a cell structure as a quotient of $T$ by a map that is affine on cells, and the image of each $1$-cell of $T$ is either a saddle connection or a point.\end{proposition}

\noindent Proof: The existence of a quotient cell structure is clear, but we must show that no triangle degenerates to a union of two distinct $1$-cells of equal slope, since this is the one possible affine quotient of a triangle that we do not allow. If this happened, the two saddle connections would have to be saddle connections of $q_\infty$, which means they are not vertical or horizontal. However, if this were the case, the parallel cells would either both have negative slope or both positive slope, and that would mean that some triangle of $T$ in $q_n$ for all sufficiently large $n$ had three sides whose slope had the same sign, which is impossible. Finally, we claim that this cell structure is a veering triangulation. The sequences pseudo-metrics $d_m^h$ and $d_m^v$ defining vertical and horizontal components of the distance between cone points converge. Thus for the two ends $a$, $b$ of a saddle connection $\gamma$ of $T$, and any third cone point $c$, the lack of cone points in the open rectangle associated to $T$ is equivalent to the condition that $$d_m^h(a,c) - d_m^h(a,b),d_m^h(b,c) - d_m^h(a,b),d_m^v(a,c) - d_m^v(a,b),d_m^v(b,c) - d_m^v(a,b)$$ cannot all be negative. This remains true when taking limits. $\Box$.

\begin{proposition} Let $\{(X_m,q_m)\}$ be as in \Cref{DegenCell}, with common triangulation $T$ and associated train track $\tau$. Then the train track $\tau_\infty$ dual to the quotient triangulation $T_\infty$ on $(X_\infty,q_\infty)$ is the complement of the inessential subgraph $H$ consisting of edges whose periods go to $0$. Moreover, the period of each saddle connection dual to an edge $e$ of $\tau$ converges along $(X_m,q_m)$ to $0$ if that edge is vanishing, and to the period of the saddle connection dual to the edge $e_\infty$ of $\tau_\infty$ that contains $e$ if $e$ is not a vanishing edge of $T$.\end{proposition}

\noindent Proof: This is immediate from the construction of the quotient map. $\Box$\\

\noindent The following observation will be essential for obtaining contraction:

\begin{lemma}\label{InTheInterior} Suppose $(X_m,q_m) \to (X_\infty,q_\infty)$ satisfy the hypotheses of \Cref{DegenCell} with common train track $\tau$. Suppose that every large edge of $\tau$ has no neighbors in the inessential subgraph $H$ whose edges are dual to the vanishing saddle connections in the dual triangulation $T$. Then the elements of $\mathcal{MF}(\tau)$ coming from $(X_m,q_m)$ converge to a point on the interior of the cone $\mathcal{MF}(\tau)$.\end{lemma}

\noindent Proof: It is sufficient to show that $(X_m,q_m)$ are given by convergent sequences of tangential measures whose limits satisfy the following: (1) every edge has positive length, and (2) the triangle inequality is strict for each complementary $3$-gon, which belongs in the limiting equivalence class. Condition (2) is equivalent to a sequence of rectangle decompositions for which every drawn separatrix of a simple zero of $(X_m,q_m)$ has height tending to a nonzero limit, since if $a,b,c$ are the three sides of such a $3$-gon and $\ell(s)$ denotes the sum of the legnths of the sides of $s$, $\ell(a) + \ell(b) - \ell(c)$ is twice the length of the separatrix corresponding to the vertex where $a$ meets $b$.\\

\noindent We first note that if we take the rectangle decomposition of \Cref{PropCanonicalRectangles} on $(X_m,q_m)$, then heights of all rectangles converge, and the rectangles whose heights have limit $0$ will be those whose small partners are all vanishing (which includes those with no small partners). The separatrices of height $0$ will be those corresponding to vertices where all three edges are vanishing. We will describe two operations that do not change the equivalence class of the tangential measure, starting with the collection of heights that come from the limit of the heights of the saddle connections in the triangulation $T$ along the sequence $(X_m,q_m)$.\\

\noindent We perform the following operation once: for each vertex $v$ belonging to a large edge $e$ with neighbors $a$ and $b$, decrease the heights of $a$ and $b$ by $1/4$ the minimum of the heights of $a$ and $b$ and increase. Each edge therefore loses at most half of its length, and each separatrix loses at most half its length. Each large edge gains height that is bounded below by $1/4$ the minimum of the heights of the rectangles corresponding to the small and mixed edges of the train track $\tau_\infty$.\\

\noindent Now, having performed the first operation, we perform the following process finitely many times, until every separatrix and every rectangle has positive height: For each large or mixed edge rectangle $R_e$ of height $h(e)$, decrease the height of $R(e)$ by $h(e)/2$ and increase the height of the neighbors of $e$ by $h(e)/2$ if $e$ is mixed, and by $h(e)/4$ if $e$ is large. This process can only increase the length of each separatrix, cannot cause an edge of nonzero length to have zero length. We can define a Lebesgue measure-preserving vertical flow $\{\phi_t|t \in \rr\}$ on a full measure subset of the collection of vertical unit tangent vectors of any $(X_m,q_m)$; in particular we define this flow on the subset of the unit tangent bundle consisting of vertical vectors that do not lie over a leaf of the vertical foliation that contains a cone point. By Poincar\'{e} recurrence, there is a flow line starting in each rectangle and returning to the same rectangle in the same direction, with respect to the trivialization of the tangent bundle determined by $q$. In particular, if the train route corresponding to such an infinite half-leaf starts somewhere in a rectangle corresponding to $a_0$ and moves toward a small end of the edge $a_0$ in such a way that it will eventually traverse $a_0$ in the same direction, the sequence of widths of rectangles traversed is initially increasing but returns to the width of $a_0$. Let $a_0,a_1,...,a_k$ be the train track edges corresponding to the longest initial sequence of rectangles of increasing width for this train route. Then $a_k$ is large, since the following edge $a_{k+1}$ is narrower (by the switch conditions and the fact that no rectangle has width zero). This means that after $k$ repetitions of the operation described, the length of $a_0$ and the lenth of the separatrix corresponding to the vertex where $a_0$ meets $a_1$ both increase. Since such a $k$ exists for every choice of small initial half edge $(a_0 \to a_1)$ it follows that finitely many interations of this procedure described gives positive height to every edge and every separatrix, as desired. $\Box$\\

\begin{proposition}\label{pseudofilling} Let $\{\tau_m\}_{m = -\infty}^\infty$ be a periodic splitting sequence of labelled train tracks, i.e. assume there is some $k$ such that $\tau_{m + k}$ is isotopic to the image of of $\tau_m$ as a labelled train track by a fixed homeomorphism $\phi$ of the underlying surface $S$ that does not depend on $m$, and $\tau_{m + 1}$ is obtained from $\tau_m$ by performing a left or right split on an edge of $\tau_m$ whose label only depends on $m$ mod $k$. Then the following are equivalent: \begin{enumerate}
    \item $\phi$ is pseudo-Anosov.
    \item The collection of edges that eventually split form a filling subtrack of $\tau_m$ for some $m$.
    \item The collection of edges that eventually split form a filling subtrack of $\tau_m$ for every $m$.
\end{enumerate} Under these equivalent conditions, we say that the sequence has \emph{pseudo-Anosov period}, and the invariant foliation of the pseudo-Anosov is supported on the collection of edges of $\tau_m$ that eventually split.\end{proposition}

\noindent Proof: First, we claim that the last two conditions are equivalent, since left and right splits preserve the genus and number of punctures of each connected component. The complement of the subgraph consisting of the edges that never split is therefore filling for some $m$ if and only if it is filling for all $m$.\\

\noindent Second, we note that if the collection of edges that split is not filling, then there is some essential curve represented by a train route intersecting only the edges that never split. This curve is invariant up to isotopy under $\phi$, so $\phi$ is not pseudo-Anosov.\\

\noindent Finally, suppose that the collection of edges that split is filling. Every train route can be continued indefinitely so it must eventually traverse an edge in the same direction twice, so there is at least some transverse measure on each $\tau_m$. We also know that the projectivization of the cones of transverse measures on the tracks $\tau_m$ have nonempty intersection by the finite intersection property of compact sets. In fact, the set of such measures is a convex cone, so by the Brouwer fixed point theorem, there is a $\phi$-invariant projective class of measured foliation on some $\tau_m$. Let $\mu_m$ be transverse measures such that $(\tau_m,\mu_m)$ are all equivalent and projectively $\phi$-invariant.\\

\noindent It is possible that the support of $\mu_m$ is a proper subset of the edge set $E(\tau_m)$, but we first claim that the support of $\mu_m$ is filling. Indeed, if it is not, let $H_m$ be the subgraph of $\tau_m$ whose edges are those for which $\mu_m$ is zero; the number of edges in $H_m$ is nondecreasing in $m$, and therefore constant by projective invariance. Thus an edge of $\tau_m$ outside the support of $\mu_m$ cannot win a split, unless both losers are both outside the support of $\mu_m$. Now, there is a maximal (up to homotopy), possibly disconnected subsurface $Y$ filled by the support of $\mu_m$, namely the union of all open disks and once-punctured disks with boundary in the support of $\mu_m$. Some boundary component of $Y$ is an essential curve $\alpha$ which can be realized by a curve transverse to $\tau_m$ intersecting only edges outside the support of $\mu_m$; the order in which these intersections occur is eventually invariant for the sequence $\{\tau_m\}$. But then, each of these edges in $H_m$ must be small at the end where it meets the support of $\mu_m$, by the switch condition and non-negativity of edge weights. This never changes under a split, so the edges crossed by this representative of $\alpha$ never split. This contradicts the assumption that the collection of edges that split is filling.\\

\noindent Now, we will go one step further, and claim that the support of $\mu_m$ is exactly the collection of edges that eventually split. Again, we let $H_m$ be the graph whose edges are the edges of $\tau_m$ that don't eventually split. Every train route of $H_m$ is finite, because no curve is carried by $H_m$. there are only finitely many train routes in $H_m$, and if an edge of $H_m$ splits, the number of train routes that only traverse edges in $H_m$ decreases. To any train route in $H_{m+1}$ we can find a train route in $H_m$ consisting of either the same sequence of edges in $H_m$, or of a train route that inserts the splitting edge $e$ between the winner and loser at opposite ends of $e$ if they appear consecutively; this is an injective map of sets of train routes, and the train routes in $H_m$ that passed through the two losers on either side of $e$ are not in the image of this injection. Thus no edge of $H_m$ splits.\\

\noindent Now, we are ready to construct our pseudo-Anosov homeomorphism. There is an induced splitting sequence of period $\leq k$ on the sequence of complements $\{\sigma_m\}$ of the graphs $\{H_m\}$; recall $\sigma_m$ is filling but possibly not complete. If there is a veering triangulation dual to these complements, we can determine the signs of the slopes of the saddle connections dual to the edges based on whether each edge most recently underwent a left split or a right split. We are not yet asserting the existence of such a triangulation, but we are free to label each edge of each $\sigma_m$ that eventually splits with the additional data ``positively sloped saddle connection" or ``negatively sloped saddle connection". (Note that for any small or mixed edge, this can be deduced from the counterclockwise order of the two vertices at a small end.) At each vertex, the two small half-edges carry opposite sign slope labels, so in particular not all edges labelled have the same sign. We may build another train track $\sigma_m^\prime$, which is \emph{isotopic as a marked graph, but not smoothly isotopic} to $\sigma_m$: at each vertex, $\sigma_m$ and $\sigma_m^\prime$ have different large half-edges labelled by signs of the same slope. The labels of the slope signs satisfy the condition needed to form a veering triangulation, but we need a system of heights and widths for all saddle connections to build a veering triangulation dual to $\sigma_m$ and $\sigma_m^\prime$ (dual as ribbon graphs). We already have widths up to scale for the saddle connections from the transverse measures $\mu_m$, so we just need to find heights.\\

\noindent The key observation is that splittings of $\sigma_m$ are dual to splittings of $\sigma_{m+1}^\prime$ as $m$ moves in the reverse direction: If $e$ is the label of an edge that is large in $\sigma_m$ and small in $\sigma_{m+1}$, then $e$ is small in $\sigma_m^\prime$ and large in $\sigma_{m+1}^\prime$. By our previous argument, the sequence $\{\sigma_m\}$ admits a sequence of compatible transverse measures of full support that are projectively invariant under $\phi^{-1}$. We now claim that we can define a half-translation structure equipped with a triangulation by saddle connections dual to $\sigma_m$, (and $\sigma_m^\prime$), by picking heights and widths of all saddle connections according to the pair of projectively invariant transverse measures on $\sigma_m$ and $\sigma_m^\prime$. Now, each split simply consists of replacing the diagonal of one convex quadrilateral whose sides have alternating slope signs with the other diagonal; it follows that after $k^\prime$ switches, we have reached a triangulation isotopic to the translation of our original triangulation by $\phi$. The heights and widths of the triangles differ by scalars $\lambda_h$ and $\lambda_w$ which are the same for all triangles. Since the area is the same, $\lambda_w\lambda_h = 1$; since the transverse measures always decrease when edges split we have $\lambda_w \neq 1$. It follows that $\phi$ is homotopic to a map with derivative $g_{\lambda_w}$ with respect to the $q$-metric for some half-translation surface structure $(X,q)$, i.e. $\phi$ is homotopic to a pseudo-Anosov map.$\Box$

\begin{definition}
Let $\tau$ be a complete filling birecurrent train track carrying a measured foliation $\mu$, and let $H$ be an inessential subgraph. Assume that for all tangential measures $\nu$ the associated quadratic differential has no vertical saddle connections. We say that $(\tau^\prime,H^\prime,\mu^\prime)$ is an \emph{improvement} of $(\tau,H,\mu)$, or \emph{improves} $(\tau,H,\mu)$, if there is a positive integer $k$ and a sequence of pairs of train tracks and inessential subtracks $(\tau,H) = (\tau_0,H_0), (\tau_1,H_1),...,(\tau_k,H_k) = (\tau^\prime,H^\prime)$ such that $\tau^\prime$ carries $\mu$, the tracks $\tau_i - H_i$ are all isotopic, and for $1 \leq i \leq k$, $(\tau_i,H_i)$ is obtained from $(\tau_{i - 1},H_{i - 1})$ by splitting a single edge of $\tau_{i - 1}$ such that at least one loser is in $H_{i - 1}$. We say $(\tau,H,\mu)$ is \emph{resolved} if it has no improvements.\end{definition}

\begin{lemma}For any complete filling birecurrent train track $\tau$ on $S_{g,n}$, inessential subgraph $H$, and transverse measure $\mu$ on $\tau$, the number of improvements of $(H,\tau,\mu)$ is bounded by a constant only depending on $g$ and $n$. In particular, every $(H,\tau,\mu)$ has a resolution.\end{lemma}

\noindent Proof: First, we show that the number of splits in which the $H_i$ does not contain the splitting edge is uniformly bounded. Then we will show that the number of improvements in which all splitting edges are in the $H_i$ is uniformly bounded.\\

\noindent Take a pole-free-disk cover of $S_{g,n}$. For each lift of a complementary region of the \emph{complement} of $H_i$, which we will call $R$, we may sum, over all sides of $R$, the number of train routes that begin with a lift of an edge of $H_i$ that is contained in $R$, whose second edge is in $\partial R$, and whose remaining edges are on the same side of $R$. This number can never increase via an improvement, and there must be some $R$ for which this decreases if the splitting edge is not in $H_i$. Since this number is finite and uniformly bounded for each $R$, and there are only finitely many such $R$ up to the action of the deck group, this proves that the number of times the splitting edge does not belong to $H_i$ is uniformly bounded. Now, the number of times that an improvement occurs entirely by splitting edges in $H_i$ is bounded because, as shown in the proof of \Cref{pseudofilling}, the number of train routes contained entirely in $H_i$ is uniformly bounded and decreases upon such splits. $\Box$.

\begin{proposition}Let $\tau$ be a train track arising from veering triangulation of a suited quadratic differential; let $\mu$ be the transverse measure on $\tau$ associated to this quadratic differential. Suppose that $(\tau,\mu)$ splits to $(\tau_*,\mu_*)$, and $\tau_*$ is, as a labelled train track, the translate of $\tau$ by a homeomorphism $\phi$ in a pseudo-Anosov mapping class. For each inessential subgraph $H$ of $\tau$ write $H_*$ for the image of $H$ in the induced isomorphism from $tau$ to $\tau_*$. Moreover, suppose that for each inessential $H$, the $\phi$-translates of improvements of $(\tau,H,\mu)$ are exactly the improvements of $(\tau_*,H_*,\mu_*)$ (up to isotopy). Let $H$ be the inessential subgraph of $\tau$ consisting of edges that never split along the sequence from $\tau$ to $\tau^\prime$ Then $(\tau,H,\mu)$ has a resolution $(\tau^\prime,\mu^\prime,H^\prime)$ with the following properties: \begin{enumerate}
    \item $H$ and $H^\prime$ have the same number of edges, and they have the same labels.
    \item If an edge that is not labelled by $H$ splits in the improvement, neither winner is labelled by $H$.
    \item $H^\prime$ does not contain any neighbors of large edges (and therefore $H^\prime$ also contains no large edges).
\end{enumerate}\end{proposition}

\noindent Proof: First, we note that every edge of the complement of $H$ will eventually split, and edges labelled by $H$ are never winners of such splits, since an edge can only return to its initial state (small, mixed or large) if it is a winner exactly twice as often as it splits. Now, we note that if we start with $(\tau,H,\mu)$ the following two types of sequence commute: (i) splitting sequences that translates by $\phi$, and (ii) improvement by splitting a single edge. It follows that any improvement $(\hat{\tau},\hat{H},\hat{\mu})$ of $(\tau,H,\mu)$ and the corresponding improvement of $(\tau_*,H_*,\mu_*)$ satisfy the same hypotheses as $(\tau,H,\mu)$ and $(\tau_*,H_*,\mu_*)$, and the edges labelled by $H$ are exactly the edges that never split as $\hat{\tau}$ splits to $\phi(\hat{\tau})$. By induction we may upgrade operation (ii) to any improvement, and it will commute with a sequence with property (i). Moreover, the edges that never split in the sequence with property (i) are exactly those labelled by $H$. So now, if we pick a resolution $(H^\prime,\tau^\prime,\mu^\prime)$, it satisfies the same hypotheses as $(\tau,H,\mu)$. Now, consider any split of an edge of $\tau^\prime$ compatible with $\mu^\prime$. The loser of that split cannot be in $H$ because $(\tau^\prime,H^\prime,\mu^\prime)$ has no improvements. The winner of the split cannot be in $H^\prime$ because that edge will eventually split in a sequence with property (i), which implies that any neighbor in $H^\prime$ must be a loser, and each edge can only split in one way that is compatible with the transverse measure. Therefore, no large edge of $\tau^\prime$ has a neighbor in $H^\prime$. $\Box$

\section{Exponential Contraction in the Hilbert Metric}

\noindent In this section, we prove the following key lemma:

\begin{lemma}\label{HilbertContraction} Let $\{g_t(X,q): 0 \leq t \leq T\}$ be a geodesic arc in $\mathcal{M}_{g,n}$ starting and ending at points with well-defined $L^\infty$ Delaunay triangulations. Let $K$ be a compact subset of $\mathcal{M}_{g,n}$ and let $\theta \in (0,1]$. Let $A_0$ and $A_T$ be the $L^\infty$ Delaunay triangulations of $A_0$ and $A_T$ respectively, with corresponding vertical tracks $\tau_0$ and $\tau_T$. Let $\lambda$ be Lebesgue measure on $\rr$. Assume $$\lambda \{t \in [0,T]: \pi(g_t(X,q)) \in K\} > \theta T.$$

Then there is a veering triangulation $A_T^\prime$, with vertical train track $\tau_T^\prime$ obtained by performing a uniformly bounded number splits of the vertical track $\tau_T$ of $A_T$ and which carries the vertical foliation of $(X,q)$, with the property that the the diameter of $\mathcal{MF}(\tau_0)$ is most $Ce^{-\alpha t}$ in the Hilbert metric on $\mathcal{MF}(\tau_T^\prime)$. The positive constants $C,\alpha$ depend only on $K$ and $\theta$.
\end{lemma}

\noindent The following proposition, when combined with \Cref{pseudofilling}, shows that splitting sequences of vertical tracks of $L^\infty$ Delaunay triangulations of recurrent Teichm\"uller geodesic segments contain pseudo-Anosov words:\\

\begin{proposition}\label{ThicketySplit}Suppose that the splitting sequence of vertical tracks of $L^\infty$ Delaunay triangulations from $(X,q)$ to $g_t(X,q)$ is such that both ends have the same labelled train tracks dual to the Delaunay triangulation, and both lie in $\mathcal{M}_{g,n}^{(\epsilon)}$. Let $\pi$ be the projection from $QD^1(\mathcal{T}_{g,n})$ to $\mathcal{M}_{g,n}$ and let $\lambda$ be Lebesgue measure on $\rr$. Then there is a constant $M(\epsilon)$ such that if the collection of edges that split along the sequence train tracks dual to the $L^\infty$ Delaunay triangulations is not filling, then $\lambda \{s \in [0,t]: \pi(g_s(X,q)) \in \mathcal{M}_{g,n}^{(\epsilon)}\} \leq M(\epsilon)$.\end{proposition}

\noindent Proof: If the collection of edges that eventually split is not filling, then there is a curve $\alpha$ intersecting only those edges whose length in the $q$-metric and the in $g_tq$ metrics are both uniformly bounded, since it can be drawn intersecting a fixed sequence of $L^\infty$ Delaunay triangles, and the diameters of these triangles are uniformly bounded on unit area quadratic differentials over $\mathcal{M}_{g,n}^{(\epsilon)}$. So there is a number $A_\epsilon$ such that $\alpha$ has transverse length at most $A_\epsilon$ with respect to the vertical and horizontal measured foliations of $(X,q)$ and $g_t(X,q)$. It follows that the horizontal length of $\alpha$ at $g_s(X,q)$ is at most $e^{s-t}A_\epsilon$ and the vertical length of $\alpha$ at $g_s(X,q)$ is at most $e^{-s}A_\epsilon$, and therefore the length of $\alpha$ in the $g_s q$-metric is at most $(e^{-s} + e^{s-t})A_\epsilon$ for $0 \leq s \leq t$. There is a constant $B_\epsilon$ such that no unit area quadratic differential lies over $\mathcal{M}_{g,n}^{(\epsilon)}$ if it has a curve of flat length less than $B_\epsilon$. However, $\alpha$ has length less than $B_\epsilon$ in the $g_s q$-metric unless $e^{-s}A_\epsilon \geq B_\epsilon/2$ or $e^{t-s}A_\epsilon > B_\epsilon/2$. Thus, unless $s < \log(2 A_\epsilon/B_\epsilon)$ or $t-s < \log(2 A_\epsilon /B_\epsilon)$, $\alpha$ is too short in the $g_s q$-metric for $g_s q$ to lie over $\mathcal{M}_{g,n}^{(\epsilon)}$. In particular, we can take $$M(\epsilon) = 2\log(2 A_\epsilon /B_\epsilon). \Box$$

\noindent Proof of \Cref{HilbertContraction}: We pick points $g_{t_k}$ along the geodesic that lie in $QD(\mathcal{M}_{g,n}^{(\epsilon)})$ whose vertical train tracks are isomorphic as labelled train tracks, but spaced such that between any two, enough of the geodesic lies in $\mathcal{M}_{g,n}^{(\epsilon)}$ the conclusion of \Cref{ThicketySplit} does not hold, i.e. for any two chosen points $g_{t_k}(X,q)$ and $g_{t_{k+1}}(X,q)$ we have $$\lambda \{s \in [t_k,t_{k+1}]: \pi(g_s(X,q)) \in \mathcal{M}_{g,n}^{(\epsilon)}\} > M(\epsilon).$$ Moreover, we assume that for all of the $g_{t_k}$ the resolutions of all inessential subgraphs are the same. By the pigeonhole principle, there is a number $c > 0$ such that we can pick at least $cT$ such points, for all sufficiently large $T$.\\

\noindent For any two such points, the corresponding train tracks differ by a pseudo-Anosov mapping class, by \Cref{pseudofilling}.

The collection of surfaces in $S \in QD^1(\mathcal{T}_{g,n}^{(\epsilon)})$ with a fixed marked Delaunay triangulation is compact. By proper discontinuity of the mapping class group action, there are only finitely many splitting sequences that move $S$ a distance at most $d$ in Teichm\"uller space, hence only finitely many labelled splitting sequences of $L^\infty$ Delaunay triangulations arising from such $L^\infty$ Delaunay triangulations starting and ending at points in $QD^1(\mathcal{T}_{g,n}^{(\epsilon)})$ distance at most $d$ apart, for each $d > 0$, by \Cref{UniquelySplit}.

By construction, the average value of $t_{k+1} - t_k$ is bounded. Thus, by the pigeonhole principle, there is some number $c^\prime$ such that the same word (the same sequence of left and right splits) appears in the splitting sequence at least $c^\prime T$ times. Now, if this pseudo-Anosov is not supported on the principal stratum, we may perform the splits in order to resolve the inessential subgraphs corresponding to non-generic singularities on the quadratic differentials on its axis. With respect to the train tracks coming from this resolution, the pseudo-Anosov meets the hypotheses of \Cref{InTheInterior} and \Cref{CoroGoodWordPrep}, and therefore meets the hypotheses of \Cref{CoroGoodWord}. Any word meeting the criteria of \Cref{CoroGoodWord} induces a contraction in the Hilbert metric on the space of equivalence classes of tangential measures. It follows that our original splitting sequence, plus perhaps a bounded number of additional splits in order to resolve an inessential subgraph the final time, contains one of finitely many words meeting the hypotheses of \Cref{CoroGoodWord} at least $c^{\prime\prime}T$ times, for some smaller constant $0 < c^{\prime\prime} < c^\prime < c$ that does not depend on $T$, and therefore its image has exponentially small diameter. $\Box$\\

\noindent Remark: There are many pseudo-Anosov words that have the same combinatorial action on $\mathcal{MF}(\tau)$, in the sense that they map extremal points of the cone into the same faces of the cone. A collection of such words, interspersed with ``noise" in the same sense as in the proof, would have similar properties to the words in \cref{CoroGoodWord}. Any lower bound for the constant $\alpha$ that could be obtained from chasing the proof will be extremely bad compared to the actual best possible value.\\

\section{Comparing Hilbert and Euclidean Metrics}

\noindent In appendix A of \cite{lhc} it was shown that the decomposition of $QD(\mathcal{M}_{g,n})$ into sets based on the $L^\infty$ Delaunay triangulation is a locally finite decomposition into convex sets (convex with respect to systems of period coordinates). For each such convex set, one may use the periods of all saddle-connections that are 1-cells of the $L^\infty$ Delaunay triangulation to get an embedding into $\cc^d$ for some $d$; we call such sets \emph{Delaunay Coordinate Charts}.

\begin{definition} We endow $QD(\mathcal{M}_{g,n})$ with the path metric induced by this collection of convex sets, i.e. the largest metric which is less than or equal to the Euclidean metric on each convex period coordinate chart. We will call this metric the \emph{Euclidean metric} and denote it by $d_E$.\end{definition}

\noindent We remark that this particular choice  of charts is arbitrary; any locally finite system of linear metrics coming from systems of periods of saddle connections gives rise to another metric in the same local Lipschitz class. This metric is not proper, but each point is contained in a ball of positive radius with compact closure.

\begin{proposition}\label{UnitLengthControl}Fix a compact subset $K$ of $\mathcal{M}_{g,n}$ and a number $N$. Let $(X,q)$ be a unit area half-translation surface with $X \in \mathcal{T}_{g,n}$ lying over $K$. We will write $e \in T$ if $e$ is an edge a triangulation in $T$ and $R(e)$ and $I(e)$ for the real and imaginary parts of the period of $e$. For any triangulations $T_1,T_2$ graph distance at most $N$ away from the $L^\infty$ Delaunay triangulation of $(X,q)$ in $\mathrm{VT}(X,q)$, the following quantities are uniformly bounded (i.e. by bounds depending only on $K$ and $N$):\begin{itemize}
    \item $\sum\limits_{e \in T_i} |R(e)|$ and $\sum\limits_{e \in T_i} |I(e)|, i = 1,2$
    \item $\left[\sum\limits_{e \in T_i} |R(e)|\right]^{-1}$ and $\left[\sum\limits_{e \in T_i} |I(e)|\right]^{-1}, i = 1,2$
    \item The entries of any period coordinate transition matrix between a system of period coordinates coming edges of $T_1$ and a system coming from edges of $T_2$.
\end{itemize}
\end{proposition}

\noindent Proof: By local finiteness of $L^\infty$ Delaunay triangulations it suffices to prove these statements just for the $L^\infty$ Delaunay triangulation itself. Clearly the lengths of the edges of $L^\infty$ Delaunay triangles are all bounded above by twice the diameter of $X$ in the $L^\infty$  $q$-metric, since the triangles are inscribed in squares $S_i$ whose centers are $L^\infty$ distance half the side length of $S_i$ to the nearest cone point (in the $L^\infty$ flat metric). This takes care of the first item on the list. Now, given a triangle in $\rr^2$ whose edges $e_1,e_2,e_3$ have widths $w_1,w_2 < w_3$ and heights $h_1 > h_2,h_3$, respectively the area of this triangle is given by

$$w_3 h_1 - [w_1h_1 + w_2h_2 + w_3h_3]/2.$$

\noindent (This just comes from drawing drawing a circumscribed oriented rectangle and subtracting off three triangular regions in the rectangle but outside the triangle.) The second item follows from the first and the fact that if the widths are all small and the heights remain bounded, the area of each triangle is small and the number of triangles remains fixed, so the area of all of the triangles cannot sum to 1.$\Box$\\

\begin{definition}If $\alpha$ is a curve on a Riemann surface $X$ with quadratic differential $q = (dx + i dy)^2$ in local coordinates, the \emph{height} of $\alpha$ on $(X,q)$ is the infimum of $\int_\gamma |dy|$ over all simple closed curves $\gamma$ representing $\alpha$. Similarly, the \emph{width} of $\alpha$ is the infimum of $\int_\gamma |dx|.$\end{definition}

\begin{proposition}\label{NormalizingFunctional}Suppose $\tau$ is obtained from the $L^\infty$ Delaunay triangulation of $(X,q)$ by performing finitely many splits compatible with the vertical foliation of $(X,q)$. Let $VC(\tau)$ be the set of vertex curves of $\tau$, which is obtained by performing at most $N$ splits of the vertical train track of $(X,q)$. Assume $(X,q) \in K$ are as in \Cref{UnitLengthControl}. Then the following are bounded (i.e. the bound depends only on $K$):\begin{itemize}
    \item $\sum\limits_{\alpha \in VC(\tau)} \mathrm{height}(\alpha)$
    \item $\left[\sum\limits_{\alpha \in VC(\tau)} \mathrm{height}(\alpha)\right]^{-1}$
\end{itemize}
\end{proposition}

\noindent Proof: For each vertex curve, there is a realization transverse to the horizontal foliation that passes through only the rectangles corresponding to the edges of the train track it carries. The parallel transport holonomy along such a curve is $\pm 1$. We describe two families of curves in the orienting double cover $\tilde{X}$: For each vertex curve $\alpha$ on which the $q$-holonomy is $-1$, concatenate this curve with itself, and take a lift to $\tilde{X}$. For each vertex curve $\alpha$ for which the $q$-holonomy is $1$, take the two distinct lifts of a representative of $\alpha$, but such that their projections to $X$ traverse the same curve in opposite directions. Each of these cycles can be picked transverse to the horizontal foliation of $\tilde{X}$ so that its tangent vector forms an acute (or $0$) angle with the upward direction at all points. The homology classes of these cycles are anti-invariant under the involution, so their Poincar\'{e} duals belong to $H_{odd}^1(\tilde{X};\rr)$. On the other hand, the transverse measure of the vertical foliation on $\tilde{X}$ is a non-negative linear combination of Poincar\'{e} duals to these vertex cycles, so they span the real parts of the local parameter space for $QD(\mathcal{T}_{g,n})$, i.e. they span $H_{odd}^1(\tilde{X};\rr).$\\

\noindent It follows that the heights of saddle connections in the $L^\infty$ Delaunay triangulation are rational linear combinations of the heights of vertex curves, and vice versa; the numerators and denominators of the transition matrices for such pairs of bases are bounded by a constant that depends only on $N$. This now follows from the upper and lower bounds on the sum of the heights of the saddle connections from \Cref{UnitLengthControl}.$\Box$

\begin{proposition}\label{eucleqhilb}Assume all notation is as in \Cref{UnitLengthControl}. Suppose that $(X,q)$ is a unit area half-translation surface lying over $K$. Let $\sigma$ be the vertical train track of $(X,q)$, and assume it splits at most $N$ times to a track $\tau$ carrying the vertical foliation of $(X,q)$. Let $\nu$ be the horizontal foliation of $(X,q)$ and let $\nu^\prime$ be another horizontal foliation realized as a tangential measure on $\tau$. Let $d_H^\mathrm{\tau,t}$ be the Hilbert distance on the space of \emph{tangential measures on} $\tau$. Let $(X^\prime,q^\prime)$ be the unit area quadratic differential on the strongly stable leaf of $(X,q)$ and on the unstable leaf determined by the projective equivalence class of $\nu^\prime$.

There are constants $C_K,\epsilon_K$, depending only on $K$ and $N$, such that the following holds: $\mathrm{if} ~ d_H^\mathrm{\tau,t}(\nu,\nu^\prime) < \epsilon_K, ~ \mathrm{then}$

$$d_E((X,q)(X^\prime,q^\prime)) < C_K d_H^\mathrm{\tau,t}(\nu,\nu^\prime).$$\end{proposition}

\noindent Proof: One can transform $(X,q)$ to $(X^\prime,q^\prime)$ in three steps: first, by changing the imaginary parts of period coordinates, second by rescaling to unit area, and third, by applying the Teichm\"uller flow. We show that each of these transformations is a motion of at most $C_K d_H^\mathrm{\tau,t}(\nu,\nu^\prime)$ in the metric $d_E$.\\

\noindent First, we note that the Hilbert metric on $[0,1]$ is less than or equal to a multiple of the Euclidean metric. Second, we note that the $q$-metric can be scaled by a positive real constant c, bounded above and away from $0$ by bounds depending only on $K$, so that the sum of the heights of the vertex curves of $\tau$ on $(X,c^2q)$ is $1$. The sum of the heights of the vertex curves is a positive functional on the cone of heights, and so a segment of Hilbert length $s$ in this cone is has length at most a bounded multiple of $s$ in the Euclidean metric $\rr^{E(\tau)}$; hence the same is true of the length of the corresponding segment in $\rr^{E(\sigma)}$. Since the area of a quadratic differential is a bilinear function of the heights and widths of $L^\infty$ Delaunay Triangulation saddle connections, it follows that if $s$ is sufficiently small, a perturbation of size $cs$ has effect $O(s)$ on the area. Thus the rescaling and Teichm\"uller flow are also $O(s)$.$\Box$

\begin{proposition}\label{DefBall}Let $K$ be as above. There exists a compact set $K^\prime$ and positive real numbers $\epsilon_0 = \epsilon_0(K) > 0, t(K) > 0$ such that the following holds: given any unit area quadratic differential $(X,q) \in K$ with a (not necessarily unique) Delaunay triantulation $\tau$, and $\epsilon < \epsilon_0$, there is an $\epsilon$-ball in the metric $d_E$ restricted to the strongly stable leaf of $g_s(X,q),$ $|s| < t$ contained entirely in $K^\prime$ and on which $\tau$ is the $L^\infty$ Delaunay triangulation.$\Box$\end{proposition}

\noindent Proof: We may pick $K^\prime$ to be any compact set with $K$ in its interior. The set on which the $L^\infty$ Delaunay triangulation is constant is locally one of finitely many convex polytopes in $\cc^d$, and the Euclidean metric on $QD(\tilde{T}_{g,n})$ is locally bi-Lipschitz to the standard metric on $\cc^d$ on each of these polytopes. Thus, by compactness, there is some $\epsilon_1(K)$ such that any ball of radius $\epsilon_1$ in the Euclidean metric about a unit area quadratic differential is contained in $K^\prime$, and moreover, the part of this ball on which $\tau$ is the $L^\infty$ Delaunay triangulation contains a ball of some radius $\epsilon_2(K,K^\prime)$. By compactness, such a ball must contain a ball of some radius $\epsilon_3(K, K^\prime) > 0$ in the Euclidean metric restricted to some strongly stable leaf. The distance from this ball to the strongly stable leaf of $(X,q)$ is bounded, again by compactness. $\Box$

\begin{theorem}\label{final} Let $\{g_t(X,q): 0 \leq t \leq T\}$ be a geodesic arc in $\mathcal{M}_{g,n}$. Let $K$ be a compact subset of $\mathcal{M}_{g,n}$ containing $(X,q)$ and $g_T(X,q)$ and let $\theta \in (0,1]$. Let $\lambda$ be Lebesgue measure on $\rr$. Assume $$\lambda \{t \in [0,T]: \pi(g_t(X,q)) \in K\} > \theta T.$$ Then there exist positive constants $C, A, \epsilon_0$ depending only on $K$, such that if $(X,q),(X_1,q_1),(X_2,q_2)$ lie on the same strongly stable leaf, and $$d_E((X_1,q_1),(X,q)),d_E((X_2,q_2),(X,q)) < \epsilon_0, ~ \mathrm{then}$$ $$\frac{d_E(g_T(X_1,q_1),g_T(X_2,q_2))}{d_E((X_1,q_1),(X_2,q_2))} < Ce^{-AT}.$$\end{theorem}

\noindent Proof: We may assume $(X,q)$ and $g_T(X,q)$ are suited and have unique $L^\infty$ Delaunay triangulations and since the general case follows by continuity.\\

\noindent The main step is to show that the derivative of $g_T$, restricted to the tangent space of the strongly stable leaf, is exponentially small with respect to a metric that is locally uniformly bi-Lipschitz to the Euclidean metric. By \Cref{HilbertContraction},  there are triangulations $B_0$ and $B_T$, uniformly bounded $VT(X,q)$-distance away from the $L^\infty$ Delaunay triangulations of $(X,q)$ and $g_T(X,q)$, with associated vertical tracks $\tau_0$ and $\tau_T$ respectively, such that a splitting sequence induces a map $$F:\mathcal{MF}(\tau_0) \to \mathcal{MF}(\tau_T)$$ satisfying the Hilbert metric $$d_H(F(x),F(y)) < C^\prime e^{-AT}d_H(x,y)$$ for some $C^\prime,A \in (0,\infty)$ depending only on $K$. By \Cref{UnitLengthControl}, \Cref{eucleqhilb}, and \Cref{DefBall} there is some $\epsilon_0 > 0$ such that a ball of radius $\epsilon_0$ in $d_E$ along the strongly stable leaf of $g_s(X,q)$ has $g_{T-s}$ image of diameter at most $Ce^{-aT}$ where $|s|$ is uniformly bounded depending on $K$. In particular, this means that if we restrict to a system of period coordinates coming from saddle connections of the $L^\infty$ Delaunay triangulations of $(X,q)$ and $g_T(X,q)$, or any period coordinate systems obtained by a uniformly bounded integral change of basis, the unit tangent space of the strongly stable leaf of $(X,q)$, is contracted by at least $Ce^{-aT}/\epsilon_0$ in the corresponding period coordinate system, where $C$ depends only on $K$. $\Box$

\section{Two Closing Lemmas}

In this section we give two closing lemmas. We remark that similar statements essentially follow from Hamenst\"adt's work in \cite{HB}, but on with a different metric which is not easy to compare to the Euclidean metric, and which is not defined everywhere. She parametrizes Teichm\"uller flow lines by points pairs of points on the boundary of the curve complex, and uses visual distance in this boundary in place of Euclidean distance. We believe that our statement is simpler and easier to use in many contexts. We will formulate our statement in terms of affine invariant manifolds on $QD(\mathcal{T}_{g,n})$.\\

\begin{definition}An \emph{affine invariant submanifold of} $QD(\mathcal{M}_{g,n})$ is the closure of an orbit of $\psltwo$ acting on a stratum of $QD(\mathcal{M}_{g,n})$ by applying the linear transformation (up to sign) on periods of all saddle connections.\end{definition}

\begin{definition}An \emph{affine invariant submanifold of} $QD(\mathcal{T}_{g,n})$ is an irreducible component of the inverse image of an affine invariant submanifold of $QD(\mathcal{M}_{g,n})$ under the natural orbifold covering.\end{definition}

\noindent By the work of \cite{EMM}, every affine invariant manifold is locally described by real linear equations in period relations in period coordinates on its stratum.

\begin{definition}A system of \emph{train track coordinates} is a chart for an open set in $\mathcal{PMF}$ consisting of a set of edge weights on a complete filling birecurrent train track.\end{definition}

\begin{proposition}\label{loccoords} Each affine invariant submanifold of $QD(\mathcal{T}_{g,n})$ is a locally finite union of solutions of systems of real linear equations in a pair of systems of train track coordinates on $\mathcal{MF} \times \mathcal{MF}$. Moreover, for any $(X,q)$ we can take the corresponding pairs of train tracks to be a bounded number of splits away from the train tracks of the $L^\infty$ Delaunay triangulations that occur in every neighborhood of $(X,q)$.\end{proposition}

\noindent Proof: In the statement we are proving, we do not require that the train tracks be compatible with each other or with the quadratic differentials in any way. Given a train track corresponding to a veering triangulation and a transverse measure, there is a finite sequence of folds that take all transverse measures on $\tau$ to the interior of the cone of transverse measures on some train track $\sigma$. Since the union of the axes of pseudo-Anosov elements are dense in every stratum, we can simply apply a power of the axis of any pseudo-Anosov whose axis has an $L^\infty$-Delaunay triangulation with axis corresponding to $\tau$ and apply \Cref{InTheInterior}. Now, with respect to train track coordinates on the folded track $\tau^\prime$, the edges with zero weights on $\tau$ are integral linear combinations on the weights of edges of $\tau^\prime$. We can find one such $\tau^\prime$ for the vertical foliation and another for the horizontal foliation. Since two measured foliations are the vertical and horizontal foliation of a quadratic differential if and only if they satisfy an open condition (that sum of their intersection numbers with every measured foliation is positive), invariance of domain implies that each stratum is locally cut out by finitely many real linear equations in train track coordinates. The result for strata implies the result for affine invariant manifolds. $\Box$\\

\noindent We will write $d_E^M$ to denote the path metric induced by $d_E$ restricted to the collection of parametrized rectifiable paths $\gamma: [0,1] \to M$ with the following property: there is some $\epsilon > 0$ such that $\gamma(t_1)$ and $\gamma(t_2)$ belong to the same affine chart. Strictly speaking, this is not a metric on $M$, but on a blowup of $M$ resolving self-intersections, whose points are of the form $(p,A)$, where $p \in M$ and $A$ is a local irreducible component of $M$ at $p$. This blowup has a neighborhood cut out by $\rr$-linear equations in period coordinates at all points. On any compact subset of $M$, the injectivity radius within each local irreducible component is bounded below with respect to the Euclidean metric.\\

\begin{corollary}\label{locunitcoords}If $M$ is an affine invariant manifold blown up along self-intersections in $QD(\mathcal{T}_{g,n})$, then in the unit area submanifold of $M$, there are local product coordinates about any $(X_0,q_0)$, such that $(X,q)$ is parametrized $(\bar{\mu},\bar{\nu},t)$ in $\mathcal{PMF} \times \mathcal{PMF} \times \rr$, where $\bar{\mu}, \bar{\nu}$ are the projective classes of the vertical and horizontal foliations $\mu$ and $\nu$ of $(X,q)$, and the real coordinate $t$ is given by $\log(i(\mu,\nu_0))$, $\nu_0$ is the vertical foliation of $(X_0)$. \end{corollary}

\begin{theorem}\label{AISclosing} Fix $\epsilon, \theta > 0$. Fix a compact subset $K$ of an affine invariant submanifold $M$ of $QD(\mathcal{M}_{g,n})$. For a unit area quadratic differential $(X_0,q_0)\in K$ write $(X_t,q_t)$ for $g_t(X,q)$. Let $\lambda$ be Lebesgue measure. There exist $T_0 > 0, \delta >\delta_0 > 0$, depending only on $K$, such that if $T > T_0,$ and there exists $\phi \in \mathrm{Mod}(S_{g,n})$ such that $d_E^M(\phi(X,q),g_T(X,q)) < \delta$, and $\lambda\{t \in [0,T]: X_t \in M_{g,n}^{(\epsilon)}\} > \theta T$, then there exists $(X^\prime,q^\prime) \in M, T^\prime \in \rr$ with $d_E^M((X_0,q_0),(X^\prime,q^\prime) < \delta_0$ and $|T^\prime - T| < \delta_0$ such that $\phi(X^\prime,q^\prime) = g_{T^\prime}(X,q).$ Moreover, the choice of $(X^\prime,q^\prime)$ is unique up to $g_t$.\end{theorem}

\noindent Proof: Using the local product structure from \Cref{locunitcoords} in two systems of train track coordinates $(\mu,\nu,t)$ and $(\hat{\mu},\hat{\nu},\hat{t})$ from two tracks giving the local product structure based at $(X_0,q_0)$ and near $(X_T,q_T)$, we may find neighborhoods $U_0$ of $(X_0,q_0) = (\mu_0,\nu_0,t_0)$ and $U_T$ of $g_T(X_T,q_T) = (\hat{\mu}_T,\hat{\nu}_T,\hat{t}_T)$ containing a $\delta$-ball about each. We will always assume $\delta_0$ is small enough that the $2\delta_0$-balls in the Euclidean distance $d_E^M$ are all contained in a fixed compact subset $K^\prime$ of $K$, and the same is true for the $2\delta$ balls in the train track coordinates. Finally, we will also assume that $\delta$ is small enough that the $\delta$-ball about each unit area quadratic differential $(Y,\alpha)$ in lying over $K$ consists of quadratic differentials is small enough that the vertical foliations are carried by a single train track which can be obtained by a bounded number of folds from some $L^\infty$ Delaunay triangulation of $(Y,\alpha)$. More precisely, choose a collection of marked train tracks $\{\tau_i\}_{i \in I}$ such the the vertical foliation of $(Z,\beta)$ belongs to $\mathcal{MF}^\pitchfork(\tau)$ whenever if $d_E^M((Y,\alpha),(Z,\beta)) < \delta$, and for each compact subset of $QD(\mathcal{T}_{g,n})$ only finitely many $\tau_i$ are needed, and the choices of $\tau_i$ are $Mod(S_{g,n})$-invariant. (If $\delta$ is small enough, then the $\delta$-ball about any point in $K$ misses an open subset of $\mathcal{PMF}$, and therefore can be carried by a single train track; compactness implies locally finitely many choices are necessary.) These choices are just to ensure that the distance in product coordinates is locally bi-Lipschitz to $d_E^M$.\\

\noindent For each $\mathcal{PMF}$ factor, fix a train track carrying all  metric on $\mathcal{MF}^\pitchfork(\tau)$, by the standard Euclidean metric $\rr^{E(\sigma)}$ and $\rr^{E(\tau)}$. In the neighborhood $U_0$, let $(X_0,q_0),$ let $$F(\mu,\nu,t) = (\phi(\mu),\phi^{-1}(\nu),t).$$ We claim that $F$ is a contraction in the $\mathcal{PMF} \times \mathcal{PMF}$ factors in a neighborhood of $(X_0,q_0)$, with a contraction constant of $Ce^{-\alpha T}$, where $\alpha$ depends only on $\theta$ and $\epsilon$, and $C$ depends only on the value of the smallest $\epsilon$ such that $\mathcal{M}_{g,n}^{(\epsilon)}$ contains $K$. Indeed, the Teichm\"uller flow does not change the projective classes of vertical and measured foliations of quadratic differentials, so if $(Y,\alpha), (Z,\beta) \in U_0$ belong to the strongly stable leaf of $(X,q)$, with horizontal measured foliations $H_Y,H_Z$, respectively, then $F(Y,\alpha)$ and $\phi^{-1}(g_T(Y_\alpha))$ have the same strongly unstable projective measured foliation, and the same is true for $F (Z,\beta)$ and $\phi^{-1}(g_T(Z_\beta))$. The derivative of $\phi$ with respect to train track coordinates is an invertible integer matrix with bounded entries, and $g_T$ contracts by a factor of $C^{\prime}e^{-\alpha T}$, where $\alpha$ and $C$ are positive real constants depend only on $\theta$ and $\epsilon$ by \Cref{final}. $F$ contracts in the $\nu$ coordinate by a factor of $Ce^{-\alpha t}$, where $C$ is at most a bounded multiple of $C^\prime$. By a similar argument, $F$ contracts in the $\mu$ coordinate by a similar factor in the train track coordinates.\\

\noindent The contraction constant can be taken arbitrary close to $0$ if $T > T_0$ and $T_0$ is assumed to be sufficiently large, so that $F$ is guaranteed to map a $\delta_0$-neighborhood of $(X,q)$ into itself. The unique fixed point and distance estimates follow from standard properties of iterated contraction mappings. $\Box$\\

\noindent The collection of pairs of measured foliations that are realized as the vertical and horizontal foliation of a quadratic differential is open in $\mathcal{MF} \times \mathcal{MF}$. Therefore we have a local product structure, and an almost identical argument gives us the following:

\begin{theorem}\label{MGNclosing} Fix $\epsilon, \theta > 0$. Fix a compact subset $K$ of $QD^1(\mathcal{M}_{g,n})$. For a unit area quadratic differential $(X_0,q_0) \in K$ write $(X_t,q_t)$ for $g_t(X,q)$. Let $\lambda$ be Lebesgue measure. There exist $T_0 > 0, \delta >\delta_0 > 0$, depending only on $K$, such that if $T > T_0,$ and there exists $\phi \in \mathrm{Mod}(S_{g,n})$ such that $d_E^M(\phi(X,q),g_T(X,q)) < \delta$, and $\lambda\{t \in [0,T]: X_t \in M_{g,n}^{(\epsilon)}\} > \theta T$, then there exists $(X^\prime,q^\prime) \in QD^1(\mathcal{M}_{g,n}), T^\prime \in \rr$ with $d_E((X_0,q_0),(X^\prime,q^\prime) < \delta_0$ and $|T^\prime - T| < \delta_0$ such that $\phi(X^\prime,q^\prime) = g_{T^\prime}(X,q).$ Moreover, the choice of $(X^\prime,q^\prime)$ is unique up to $g_t$. $\Box$\end{theorem}

\noindent In both \Cref{AISclosing} and \Cref{MGNclosing} we conclude that $\phi$ is the pseudo-Anosov and $T^\prime$ is its dilatation. What is relevant for \Cref{AISclosing} is that the geodesic is contained in the affine invariant manifold.\\

\bibliographystyle{amsalpha}
\bibliography{pseudopigeon}

\end{document}